\documentclass[journal]{IEEEtran}
\IEEEoverridecommandlockouts                              
\let\labelindent\relax
 \usepackage{etoolbox}
\makeatletter
\patchcmd{\@begintheorem}{\textit}{\textbf}{}{}
\makeatother
 \newtheorem{definition}{\bf Definition}

  \newtheorem{Corollary}{\bf Corollary}
  \newtheorem{thm}{\bf Theorem}

  \newtheorem{remark}{\bf Remark}
 \newtheorem{lemma}{\bf Lemma}
  \newtheorem{prop}{\bf Proposition}
  \usepackage{dirtytalk}
\usepackage{amsmath}
\usepackage{amsfonts}
\usepackage{amssymb}
\usepackage{bbm}
\usepackage{framed}
\usepackage{breqn}
\usepackage{subcaption}
\usepackage{cite}
\usepackage{enumitem}
\usepackage{cases}
\usepackage{pgfplots}
\usepackage{tabu}
\usepackage{url}
\usepackage{multicol}
\usepackage{color}
\usepackage[bottom]{footmisc}
\usepackage{environ}
\usepackage{tikz}
\usetikzlibrary{arrows,decorations.markings}
\usetikzlibrary{decorations.pathreplacing}
\usetikzlibrary{backgrounds}
\usetikzlibrary{patterns}
\pgfplotsset{compat=newest}

\usepackage{nicefrac}
\usepackage{standalone}

\newcommand{\obsa}[2]{ \fill[blue, opacity=0.2] (#1+0.15,#2+0.15) rectangle (#1+1.85,#2+1.85);}

\newcommand{\obsb}[2]{ \fill[red!80!black, opacity=0.2] (#1+0.15,#2+0.15) rectangle (#1+1.85,#2+1.85);}

\newcommand{\obsc}[2]{ \fill[orange, opacity =0.2] (#1+0.15,#2+0.15) rectangle (#1+1.85,#2+1.85);}

\newcommand{\obsd}[2]{ \fill[orange!60!green, opacity=0.3] (#1+0.15,#2+0.15) rectangle (#1+1.85,#2+1.85);}
\newcommand{\obsf}[2]{ \fill[magenta!60!black, opacity=0.3] (#1+0.15,#2+0.15) rectangle (#1+1.85,#2+0.85);}
\newcommand{\obsg}[2]{ \fill[green!50!black, opacity=0.2] (#1+0.15,#2+0.15) rectangle (#1+1.85,#2+0.85);}
\newcommand{\obsq}[2]{ \fill[yellow!30!black, opacity=0.3] (#1+0.15,#2+0.15) rectangle (#1+0.85,#2+0.85);}
\newcommand{\obsqq}[2]{ \fill[blue!50!black, opacity=0.3] (#1+0.15,#2+0.15) rectangle (#1+0.85,#2+1.85);}
\newcommand{\obsqqq}[2]{ \fill[orange!60!black, opacity=0.3] (#1+0.15,#2+0.15) rectangle (#1+0.85,#2+1.85);}

\newcommand{\StaticObstacle}[2]{ \fill[red] (#1+0.1,#2+0.1) rectangle (#1+0.9,#2+0.9);}
\newcommand{\initialstate}[2]{ \fill[black!30!brown] (#1+0.1,#2+0.1) rectangle (#1+0.9,#2+0.9);}
\newcommand{\goalstate}[2]{ \fill[black!50!green] (#1+0.1,#2+0.1) rectangle (#1+0.9,#2+0.9);}

\usetikzlibrary{automata,arrows,positioning,calc}
\title{Entropy Maximization for Partially Observable Markov Decision Processes}
\author{$\text{Yagiz Savas}^{\star}$, $\text{Michael Hibbard}^{\star}$, Bo Wu, Takashi Tanaka, and Ufuk Topcu 
\thanks{$\star$ Y. Savas and M. Hibbard contributed equally to this work.}
\thanks{This work is supported in part by the grants AFRL FA9550-19-1-0169, DARPA D19AP00004, DARPA D19AP00078, FOA-AFRL-AFOSR-2019-0003, and NSF Award 1944318.}
\thanks{ All authors are with the Department of Aerospace Engineering
and Engineering Mechanics, and the Oden Institute for Computational
Engineering and Sciences, University of Texas, Austin, 201 E 24th
St, Austin, TX 78712. email: {\tt\small $\{$yagiz.savas, mwhibbard, bwu3, ttanaka, utopcu$\}$@utexas.edu}}}

\begin{document}
\maketitle

\begin{abstract} 
We study the problem of synthesizing a controller that maximizes the
entropy of a partially observable Markov decision process (POMDP)
subject to a constraint on the expected total reward. Such a controller
minimizes the predictability of an agent's trajectories to an
outside observer while guaranteeing the completion of a task expressed
by a reward function. We first prove that an agent with partial
observations can achieve an entropy  at most as well as an agent with
perfect observations. Then, focusing on
finite-state controllers (FSCs) with deterministic memory transitions,
we show that the maximum entropy of a POMDP is lower bounded by the
maximum entropy of the parametric Markov chain (pMC) induced by such
FSCs. This relationship allows us to recast the entropy maximization
problem as a so-called parameter synthesis problem for the induced pMC.
We then present an algorithm to synthesize an FSC that locally
maximizes the entropy of a POMDP over FSCs with the same number of
memory states. In numerical examples, we illustrate the relationship
between the maximum entropy, the number of memory states in the FSC,
and the expected reward.
\end{abstract}

\section{Introduction}
Entropy \cite{Cover} is an information-theoretic measure to quantify the  unpredictability of outcomes in a random variable. In this paper, we consider a sequential decision-making framework of partially observable Markov decision processes (POMDPs) in which a reward in terms of the entropy is introduced in addition to the classical state-dependent reward. More specifically, in the POMDP formulation that we consider, we look for a controller that maximizes the entropy reward while ensuring that the expected state-dependent reward is above a given threshold. Intuitively, the entropy reward plays a role to  promote the  unpredictability  of the controlled process to an outside observer.  Therefore, the considered POMDP formulation provides a meaningful framework for sequential decision-making in stochastic environments with imperfect information and nondeterministic choices, where a given task should be accomplished in the most unpredictable way. 

A controller in a POMDP resolves the nondeterminism and induces a stochastic process. 
Following \cite{Biondi,savas2018entropy}, we quantify the unpredictability of realizations in an induced stochastic process by defining the entropy of the process as the joint entropy of a sequence of random variables. We then mathematically show that the maximum entropy of a POMDP is upper bounded by the entropy of its  corresponding fully observable counterpart, which is a Markov decision process (MDP).

For a given POMDP, the main objective of this paper is to synthesize a controller that induces a process whose realizations accumulate rewards in the most unpredictable way to an outside observer. Controller synthesis problems for POMDPs are notoriously hard to solve. The optimal controllers are often required to take the full observation history into account which  makes searching for them undecidable in the infinite horizon case and PSPACE-complete in the finite horizon case \cite{undecidability,chatterjee2016decidable}. For computational tractability, POMDP controllers are often restricted to have finite states that represent finite observation memory \cite{meuleau1999solving}. Furthermore, in contrast to classical POMDP problems in which the optimal controllers are deterministic,  problems adopting information-theoretic performance criteria, such as entropy, typically admit randomized controllers that specify probability distributions over action selection. 

In this paper, we synthesize a randomized finite-state controller (FSC) for a POMDP that specifies a probability distribution over actions for each of its memory states \cite{Poupart}. In particular, we consider the POMDP entropy maximization problem over all FSCs with a fixed number of memory states. A key observation is that one can use a parametric Markov chain (pMC) to succinctly represent the product between  a POMDP and the set of all FSCs with a fixed number of memory states \cite{Cubuktepe:10.1007/978-3-030-01090-4_10, Parametric}.  By restricting our attention to FSCs with deterministic memory transitions, we recast the POMDP controller synthesis problem as a so-called parameter synthesis problem for a pMC whose entropy we aim to maximize. To build a connection between the entropy of the POMDP and that of the corresponding pMC, we first prove that the maximum entropy of a pMC induced from a POMDP by FSCs with deterministic memory transitions is a lower bound on the maximum entropy of the POMDP. Furthermore, for some specific memory transition functions in FSCs, we show that one can monotonically obtain stochastic processes induced from a POMDP with higher entropy by increasing the number of memory states in the FSCs. Finally, we present a computation algorithm, based on a nonlinear optimization problem, to synthesize parameters in an FSC to maximize the entropy of a pMC subject to expected reward constraints.

One application of this theoretical framework is the synthesis of a controller for an autonomous agent carrying out a mission in an adversarial environment. In particular, if the agent's sensor measurements are noisy and the mission is defined in terms of a reward function, the synthesized controller leaks the minimum information about the agent's trajectories to an outside observer while guaranteeing the completion of the task. Furthermore, the proposed methods can be used to distribute traffic assignments over a network with possibly noisy traffic information, which is known as stochastic traffic assignment \cite{akamatsu1996cyclic}, since higher entropy in this scenario promotes the use of different paths. 

\textbf{Related Work.} 
A preliminary version of this paper has appeared in \cite{hibbard2019unpredictable}, where we present solutions for entropy maximization over FSCs with a \textit{specific memory transition function} and the same number of memory states. This considerably extended version includes detailed proofs for all theoretical results, a nonlinear optimization problem formulating the entropy maximization over all deterministic FSCs with the same number of memory states, and an extended numerical examples section.

In a recent study \cite{savas2018entropy}, we showed that an entropy-maximizing controller for a fully observable MDP can be synthesized efficiently by solving a convex optimization problem. Moreover, we established that, when the maximum entropy of an MDP is finite, it is sufficient to focus only on memoryless controllers to induce a process with maximum entropy. It is known \cite{papadimitriou} that the synthesis of a controller that accumulates a desired level of total reward in a POMDP is, in general, intractable and such a controller typically utilizes memory. Therefore, the consideration of partial observability in the system model dramatically changes the complexity of the problem at hand. As a result, in this paper, we focus on finite-state controllers and present a nonlinear optimization problem with bilinear constraints to synthesize entropy-maximizing controllers.

In POMDPs, entropy has often been used for active sensing applications \cite{kreucher2005sensor, lauri2014stochastic,eidenberger2010active,ryan2010particle,spaan2015decision}, where an agent seeks to select actions that decrease its uncertainty on the environment by taking actions that minimize the entropy of a probability distribution. Such a distribution typically expresses the agent's belief on the task-relevant aspects of the environment. In this paper, we consider an agent that aims to maximize the entropy of its \textit{true state trajectories} instead of minimizing the entropy of its \textit{final belief distribution} on task-relevant aspects. Therefore, despite the similarity of the information-theoretic measures considered, the problem studied in this paper and the developed solution approach are significantly different from the ones investigated in active sensing literature.

In the reinforcement learning literature, the entropy of a controller has been used as a regularization term in an agent's objective to balance the trade-off between exploration and exploitation \cite{haarnoja2018composable}. As discussed in \cite{haarnoja2017reinforcement}, using a controller with high entropy, an agent can learn various ways of completing a task, leading to a greater robustness when subsequently fine-tuned to specific scenarios. The aforementioned work concerns the synthesis of a controller that balances the accumulated reward and the entropy of the induced process in a fully observable setting. Here we aim to synthesize a controller that maximizes the entropy in a partially observable setting while ensuring the accumulation of a desired level of total reward.


A range of solution techniques exist for POMDP controller synthesis using FSCs. For deterministic FSCs, existing approaches include branch-and-bound method \cite{meuleau1999solving}, automaton learning-based method \cite{zhang2015learning}, and expectation-maximization \cite{pajarinen2011periodic}. They mainly target for finding an optimal transition structure of the FSC.  As for randomized FSCs, in addition to the transition structure, one also needs to optimize the probabilistic transition probabilities between FSC states and the action selection probabilities.  To this end, researchers propose solutions using policy iteration \cite{hansen1998solving,Poupart}, gradient descent \cite{Meuleau}, and nonlinear optimization \cite{amato2010optimizing,junges2017permissive}. However, the results mentioned above only consider state-dependent reward optimization or the satisfaction of a given specification. In contrast, we consider the synthesis of FSCs for entropy maximization, which is a nonlinear objective that requires a new optimization formulation as well as solution techniques. 

\textbf{Contribution.}
The contributions of this paper are four-fold. First, we prove that the maximum entropy for a POMDP is bounded by the maximum entropy of its underlying fully observable MDP. Secondly, by restricting the scope of the FSC synthesis problem to FSCs with deterministic memory transitions, we prove that the maximum entropy of the induced pMC is a lower bound on the maximum entropy of the POMDP. Thirdly, we present a nonlinear optimization problem whose solution provides a controller that maximizes the entropy of the POMDP over all deterministic FSCs with the same number of memory states. Lastly, for deterministic FSCs, we propose a specific memory transition function which increases the entropy of the induced stochastic process with respect to an increasing number of memory states. 

\textbf{Organization.}
We provide the  modeling framework and preliminary definitions in Section \ref{section:Preliminaries}. We then formally state the entropy maximization problem for the finite and infinite horizons in Section \ref{section:Problem Statement}. We show that the maximum entropy of a POMDP is upper bounded by that of its underlying MDP in Section \ref{section:An Upper Bound on Maximum Entropy}. In Section \ref{section:Entropy Maximization Over Finite-State Controllers}, we focus on FSCs and prove that the maximum entropy of the pMC induced by a deterministic FSC is a lower bound for the maximum entropy of the POMDP. We then present a procedure to synthesize a local optimal FSC to maximize the entropy of a POMDP subject to reward constraints. We provide numerical examples in Section \ref{section:Numerical Examples} and conclude with possible future directions in Section \ref{section:Conclusions and Future Extensions}. Proofs for all technical results are provided in Appendix \ref{section:Appendix}. 

\section{Preliminaries}\label{section:Preliminaries}
We denote the power set and cardinality of a set $\mathcal{S}$ by $2^{\mathcal{S}}$ and $\lvert \mathcal{S} \rvert$, respectively. Set of all probability distributions on a finite set $\mathcal{S}$, i.e., all functions $f $$:$$ \mathcal{S}$$ \rightarrow $$[0,1]$ such that $\sum_{s\in \mathcal{S}} f(s)$$=$$1$, is denoted by $\Delta(\mathcal{S})$. 
 For a sequence $\{X_t, t$$\in$$\mathbb{N}\}$, a subsequence $(X_k, X_{k+1},\ldots, X_l)$ is denoted by $X_k^l$. The subsequence $(X_1, X_2,\ldots, X_l)$ is simply denoted by $X^l$.
 
 \subsection{Partially Observable Markov Decision Processes}
{\setlength{\parindent}{0cm}
\begin{definition}\label{def:MDP}
	A \textit{partially observable Markov decision process} (POMDP) is a tuple $\mathbf{M}=(\mathcal{S},s_I,\mathcal{A},P,\mathcal{Z},O,R)$ where $\mathcal{S}$ is a finite set of states, $s_I$$\in$$\mathcal{S}$ is a unique initial state, $\mathcal{A}$ is a finite set of actions, $P:\mathcal{S}$$\times$$\mathcal{A}$$\rightarrow$$\Delta(\mathcal{S})$ is a transition function, $\mathcal{Z}$ is a finite set of observations, $O:\mathcal{S}$$\rightarrow$$\Delta(\mathcal{Z})$ is an observation function, and $R:\mathcal{S}$$\times$$\mathcal{A}$$\rightarrow$$\mathbb{R}$ is a reward function.
\end{definition}}
For simplicity, we assume that all actions $a$$\in$$\mathcal{A}$ are available in all states $s$$\in$$\mathcal{S}$. 
For the ease of notation, we denote the transition probability $P(s'|s,a)$ and the observation probability $O(z|s)$ by $P_{s,a,s'}$ and $O_{s,z}$, respectively. 

For a POMDP $\mathbf{M}$, the \textit{corresponding fully observable} MDP $\mathbf{M}_{fo}$ is obtained by setting $\mathcal{Z}$$=$$\mathcal{S}$ and $O_{s,s}$$=$$1$ for all $s$$\in$$\mathcal{S}$. 

A \textit{system history} of length $t$$\in$$\mathbb{N}$ for a POMDP $\mathbf{M}$ is a sequence $h^t$$=$$(s_I,a_1,s_2,a_2,s_3,\ldots,s_t)$ of states and actions such that $P_{s_k,a_k,s_{k+1}}$$>$$0$ for all $k$$\in$$\mathbb{N}$. We denote the set of all system histories of length $t$ by $\mathcal{H}^{t}$. For any system history $h^t$$=$$(s_I,a_1,s_2,\ldots,s_t)$ of length $t$, there is an associated \textit{observation history} $o^t$$=$$(z_1,a_1,z_2,\ldots,z_t)$ of length $t$$\in$$\mathbb{N}$ where $O_{s_k,z_k}$$>$$0$ for all $k$$\in$$\mathbb{N}$. Note that there are, in general, multiple observation histories that are admissible for a given system history $h^t$. Finally, we denote the set of all observation histories of length $t$ by $\mathcal{O}^{t}$.
{\setlength{\parindent}{0cm}
\begin{definition}
A \textit{controller} $\pi$$:$$ \cup_{t\in \mathbb{N}}\mathcal{O}^t$$\rightarrow$$\Delta(\mathcal{A})$ is a mapping from observation histories to distributions over actions. For a POMDP $\mathbf{M}$, we denote the set of all controllers by $\Pi(\mathbf{M})$.
\end{definition}}

The probability with which the controller $\pi$ takes the action $a$$\in$$\mathcal{A}$ upon receiving the history $o^t$$\in$$\mathcal{O}^{t}$ is denoted by $\pi(a| o^t)$. 

\subsection{Entropy of Stochastic Processes}
The \textit{entropy of a random variable} $X$ with a countable support $\mathcal{X}$ and probability mass function (pmf) $p(x)$ is
\begin{align}
    H(X):=-\sum_{x\in\mathcal{X}}p(x)\log p(x).
\end{align}

We use the convention that $0$$\log$$0$$=$$0$. Let $(X_1,X_2)$ be a pair of random variables with the joint pmf $p(x_1,x_{2})$ and the support $\mathcal{X}\times \mathcal{X}$. The \textit{joint entropy} of $(X_1,X_2)$ is 
\begin{align}
\label{joint_entropy}
H(X_1,X_2):= -\sum_{x_1\in \mathcal{X}}\sum_{x_{2}\in \mathcal{X}}p(x_1,x_{2})\log p(x_1,x_{2}),
\end{align}\noindent
and the \textit{conditional entropy} of $X_2$ given $X_1$ is
\begin{align}
\label{conditional_entropy}
&H(X_2 | X_1):=-\sum_{x_1\in \mathcal{X}}\sum_{x_{2}\in \mathcal{X}}p(x_1,x_{2})\log p(x_2 |x_1).
\end{align}\noindent
The definitions of the joint and conditional entropy extend to collections of $k$$\in$$\mathbb{N}$ random variables as shown in \cite{Cover}. 
A discrete \textit{stochastic process} $\mathbb{X}$ is a discrete time-indexed sequence of random variables, i.e., $\mathbb{X}$$=$$\{X_t$$\in$$\mathcal{X}$$:$$t$$\in$$\mathbb{N}\}$.
{\setlength{\parindent}{0cm}
\noindent\begin{definition} \cite{Biondi_thesis}
The \textit{entropy of a stochastic process} $\mathbb{X}$ is defined as 
 \begin{align}\label{entropy_def_stochastic}
 H(\mathbb{X}) :=\lim_{t\rightarrow \infty}H & ( X_1,X_2,\ldots,X_t).
 \end{align}
 \end{definition}}
 The above definition is different from the \textit{entropy rate} of a stochastic process, which is defined as $\lim_{t\rightarrow \infty}\frac{1}{t}H( X^t)$ when the limit exists \cite{Cover}. The limit in \eqref{entropy_def_stochastic} either converges to a nonnegative number or diverges to positive infinity \cite{Biondi_thesis}. 

For a POMDP $\mathbf{M}$, a controller $\pi$$\in$$\Pi(\mathbf{M})$ induces a discrete stochastic process $\{S_t$$\in$$\mathcal{S}$$:$$t$$\in$$\mathbb{N}\}$ in which each $S_t$ is a random variable over the state space $\mathcal{S}$. We denote the entropy of a POMDP $\mathbf{M}$ under a controller $\pi$$\in$$\Pi(\mathbf{M})$ by $H^{\pi}(\mathbf{M})$. 

\section{Problem Statement}\label{section:Problem Statement}
We consider an \textit{agent} whose behavior is modeled as a POMDP and an \textit{outside observer} whose objective is to infer the states occupied by the agent in the future from the states occupied in the past. Being aware of the observer's objective, the agent aims to synthesize a controller that minimizes the predictability of its future states while ensuring that the expected total reward it collects exceeds a specified threshold. 

We measure the predictability of the agent's future states by the entropy of the underlying stochastic process. The rationale behind this choice can be better understood by recalling (see, e.g., Theorem 2.5.1 in \cite{Cover}) that, for an arbitrary controller $\pi$$\in$$\Pi(\mathbf{M})$, the identity
\begin{align}
    H^{\pi}(S_1, S_2,\ldots, S_N)=&H^{\pi}(S^N_t|S^{t-1})+H^{\pi}(S^{t-1})\label{derive_traj}
\end{align}
holds for any $N$$\in$$\mathbb{N}$ and $t$$\leq$$N$.
Therefore, by maximizing the value of the left hand side of \eqref{derive_traj}, one maximizes the entropy of all future sequences $(S_t,\ldots,S_N)$ for any given history of sequence $(S_1,\ldots,S_{t-1})$. 

We first consider an agent with a finite decision horizon.
\noindent\textbf{Problem 1 (Finite horizon entropy maximization):} For a POMDP $\mathbf{M}$, a finite decision horizon $N$$\in$$\mathbb{N}$, and a reward threshold $\Gamma$$\in$$\mathbb{R}$, synthesize a controller $\pi^{\star}$$\in$$\Pi(\mathbf{M})$ that solves the following problem:
\begin{subequations}
\begin{align}\label{objective1}
    &\underset{\pi\in\Pi(\mathbf{M})}{\text{maximize}}\ \ H^{\pi}(S_1,S_2,\ldots,S_N)\\ \label{constraint1}
    &\text{subject to:} \ \ \mathbb{E}^{\pi}\Big[\sum_{t=1}^{N}R(S_t,A_t)\Big]\geq \Gamma.
\end{align}
\end{subequations}

In the finite horizon entropy maximization problem, we seek a controller that randomizes the agent's \textit{finite length state trajectories} by using only the \textit{observation history} information. 

Next, we consider an agent with infinite decision horizon whose aim is to randomize its infinite length state trajectories. When the decision horizon is infinite, the total reward collected by the agent, as well as the entropy of the underlying stochastic process, may be infinite \cite{puterman2014markov,savas2018entropy}. A common approach to ensure the finiteness of the solution in infinite horizon models is to discount the collected rewards and the gained entropy in the future \cite{zhou2017infinite,mei2020global,hansen2006robust}. Accordingly, noting that
\begin{align} \label{summation_formula}
    \hspace{-0.2cm} \lim_{t\rightarrow \infty} H^{\pi}(S_1,S_2,\ldots,S_t)=&H^{\pi}(S_1)+\sum_{t=2}^{\infty}H^{\pi}(S_{t}|S^{t-1}),
\end{align}
we treat each term $H^{\pi}(S_{t}|S^{t-1})$ as a virtual entropy reward for the agent. Note that $H^{\pi}(S_1)$$=$$H^{\pi'}(S_1)$ for any $\pi,\pi'$$\in$$\Pi(\mathbf{M})$ since the initial state distribution in a POMDP is fixed. 
By discounting the agent's future rewards $R(S_t,A_t)$ as well as its virtual entropy reward $H^{\pi}(S_{t}|S^{t-1})$, we define the infinite horizon entropy maximization problem as follows: 
\noindent\textbf{Problem 2 (Infinite horizon entropy maximization):} For a POMDP $\mathbf{M}$, a discount factor $\beta$$\in$$[0,1)$, and a reward threshold $\Gamma$$\in$$\mathbb{R}$, synthesize a controller $\pi^{\star}$$\in$$\Pi(\mathbf{M})$ that solves the following problem:
\begin{subequations}
\begin{align}\label{objective2}
    &\underset{\pi\in\Pi(\mathbf{M})}{\text{maximize}}\ \ \sum_{t=2}^{\infty}\beta^{t-2}H^{\pi}(S_{t} | S^{t-1})\\ \label{constraint2}
    &\text{subject to:} \ \ \mathbb{E}^{\pi}\Big[\sum_{t=1}^{\infty}\beta^{t-1}R(S_t,A_t)\Big]\geq \Gamma.
\end{align}
\end{subequations}

For a reward function $R$$:$$S$$\times$$\mathcal{A}$$\rightarrow$$ \mathbb{R}$, which is independent of the controller $\pi$, it is known \cite{astrom} that
\begin{align*}
    \sup_{\pi \in \Pi(\mathbf{M})}\mathbb{E}^{\pi}\Big[\sum_{t=1}^{N}{R}(S_t,A_t)\Big]\leq \sup_{\pi \in \Pi(\mathbf{M}_{fo})}\mathbb{E}^{\pi}\Big[\sum_{t=1}^{N}{R}(S_t,A_t)\Big].
\end{align*}
The above inequality shows that an agent with perfect observations can collect a total reward that is at least as high as the total reward collected by an agent with imperfect observations. Since the virtual entropy rewards $H^{\pi}(S_{t} | S^{t-1})$ are functions of the agent's policy, it is not obvious whether a similar claim holds for the entropy maximization problems. In the next section, we establish that an agent with perfect observations can indeed randomize its state trajectories at least as well as an agent with imperfect observations.

It is known that deciding the existence of a controller that satisfies the constraint \eqref{constraint1} is, in general, PSPACE-complete \cite{papadimitriou}. Moreover, for $\beta$$\in$$ [0,1)$, the existence of a policy that satisfies the constraint \eqref{constraint2} is, in general, undecidable \cite{undecidability}. Therefore, the synthesis of globally optimal controllers that solve the entropy maximization problems is, in general, intractable. In the second part of the paper, we restrict our attention to a special class of controllers, namely finite-state controllers (FSCs). We present a method to synthesize FSCs that are local optimal solutions to entropy maximization problems among all FSCs with fixed number of memory states and fixed memory transition functions.

\section{An Upper Bound on Maximum Entropy}\label{section:An Upper Bound on Maximum Entropy}
In this section, we prove that an agent with perfect observations can randomize its trajectories at least as well as an agent with imperfect observations. The presented proof is based on the principle of optimality \cite{puterman2014markov} and paves the way for the formulation of an optimization problem using which we synthesize entropy-maximizing controllers.

Formally, we show that, for any $N$$\in$$\mathbb{N}\cup \{\infty\}$,
\begin{align*}
    \sup_{\pi \in \Pi(\mathbf{M})}H^{\pi}(S_1, S_2,\ldots, S_N)\leq \sup_{\pi \in \Pi(\mathbf{M}_{fo})}H^{\pi}(S_1, S_2,\ldots, S_N).
\end{align*}
We first prove that the above inequality holds for $N$$\in$$\mathbb{N}$. Then, using a monotonocity argument, we show that the inequality still holds as $N$$\rightarrow$$\infty$. 

For a given system history $h^t$$=$$(s_I,a_1,s_2,a_2,s_3,\ldots,s_t)$, let the sequences $s^t$$=$$(s_I,s_2,\ldots,s_t)$ and $a^t$$=$$(a_1,a_2,\ldots,a_t)$ be the corresponding state and action histories of length $t$, respectively. We denote the set of all state and action histories of length $t$ by $\mathcal{SH}^t$ and $\mathcal{AH}^t$, respectively. 

It can be shown that, for a POMDP $\mathbf{M}$ under a controller $\pi$$\in$$\Pi(\mathbf{M})$, the realization probability $Pr^{\pi}(s^{t+1} | s^{t})$ of the state history $s^{t+1}$$\in$$\mathcal{SH}^{t+1}$ for a given $s^t$$\in$$\mathcal{SH}^{t}$ is
\begin{align}
    Pr^{\pi}(s^{t+1} | s^{t})=\sum_{a^t\in\mathcal{AH}^t} \prod_{k=1}^t \mu_k(a_k| h^k) P_{s_t,a_t,s_{t+1}}.
\end{align}
In the above equation, $h^k$ are prefixes of $h^t$ from which the state sequence $s^t$ is obtained, and $\mu_t$ $:$ $\mathcal{H}^t$$\rightarrow$$\Delta(\mathcal{A})$ is a mapping such that
\begin{align}\label{MDP_controller}
   \mu_t(a|h^t):= \sum_{o^{t} \in \mathcal{O}^{t}} \pi(a|o^{t}) Pr(o^{t}|h^{t}).
\end{align}
We note that, for $t$$=$$1$, we have  $Pr(o^{1}$$=$$z_1|h^{1})$$=$$O_{s_I,z_1}$, and for all $t$$\geq$$2$, $Pr(o^{t}|h^{t})$ can be recursively written as
\begin{align}\label{eq:obsSeqRecursive}
    Pr(o^{t}|h^{t}) = {O}_{s_{t},z_{t}}{P}_{s_{t-1},a_{t-1},s_{t}}Pr(o^{t-1}|h^{t-1}).
\end{align}
\noindent

For a given controller $\pi$$\in$$\Pi(\mathbf{M})$ and a constant $N$$\in$$\mathbb{N}$, let ${V}_{t,N}^{\pi}$ $:$ $\mathcal{SH}^t$$\rightarrow$$\mathbb{R}$ be the \textit{value function} such that, for all $t$$<$$N$,
\begin{align}
    {V}_{t,N}^{\pi}(s^{t}) := \sum_{k=t}^{N-1} H^{\pi}(S_{k+1}|S^k_{t},S^{t}=s^{t}).  \label{eq:defEntropyMDPsum}
\end{align}
 {\setlength{\parindent}{0cm}
\noindent \begin{lemma}\label{writeValueFunction} For a POMDP ${\mathbf{M}}$, a controller $\pi$$\in$$\Pi(\mathbf{M})$, and a finite constant $N$$\in$$\mathbb{N}$, the value function ${V}_{t,N}^{\pi}$, defined in \eqref{eq:defEntropyMDPsum}, satisfies the equality  
    \begin{align}\label{eq:EntropyVfRecursive}\hspace{-0.3cm}
       {V}_{t,N}^{\pi}(s^{t}) =& H^{\pi}(S_{t+1}|S^{t}=s^{t}) \nonumber \\ 
        & + \sum_{s^{t+1} \in \mathcal{SH}^{t+1}} Pr^{\pi}(s^{t+1}| s^{t}) {V}_{t+1,N}^{\pi}(s^{t+1}) 
    \end{align}
for all $t$$<$$N$ and $s^{t} \in \mathcal{SH}^{t}$.
\end{lemma}}

We remind the reader that the proof of all technical results, including the proof of Lemma \ref{writeValueFunction}, are provided in Appendix \ref{section:Appendix}. 

For $t$$<$$N$, let ${V}^{\star}_{t,N}$ $:$ $\mathcal{SH}^t$$\rightarrow$$\mathbb{R}$ be a function such that
\begin{align}
    {V}_{t,N}^{\star}(s^{t}) := \sup_{\pi\in\Pi(\mathbf{M})}{V}_{t,N}^{\pi}(s^{t}).
\end{align}
Using Lemma \ref{writeValueFunction}, together with \textit{the principle of optimality} \cite[Chapter 4]{puterman2014markov}, we conclude that, for all $t$$<$$N$ and $s^{t} \in \mathcal{SH}^{t}$, 
    \begin{align}\hspace{-0.3cm}\label{bellman_result}
        {V}_{t,N}^{\star}(s^{t}) =&\sup_{\pi\in\Pi(\mathbf{M})} \Big[H^{\pi}(S_{t+1}|S^{t}=s^{t}) \nonumber \\ 
        & + \sum_{s^{t+1} \in \mathcal{SH}^{t+1}} Pr^{\pi}(s^{t+1}| s^{t}) {V}_{t+1,N}^{\star}(s^{t+1})\Big]. 
    \end{align}
Then, the summation in \eqref{summation_formula}, together with the definition of the value function in \eqref{eq:defEntropyMDPsum}, implies that, for any $N$$\in$$\mathbb{N}$, we have 
\begin{align*}
\sup_{\pi \in \Pi(\mathbf{M})}H^{\pi}(S_1,S_2,\ldots, S_N)={V}^{\star}_{1,N}(s_I).
\end{align*}

 The derivations presented above shows that an agent having access to state histories $s^t$ can synthesize an entropy-maximizing controller by recursively computing the values ${V}_{t,N}^{\star}(s^{t})$ via dynamic programming. In a POMDP, only observation histories are available to the agent; hence, the above derivations cannot be directly used for controller synthesis. In the next section, we consider finite-state controllers whose memory states represent statistics of the state histories and present a tractable controller synthesis method by utilizing the results of Lemma \ref{writeValueFunction}.

We now prove that the maximum entropy of a POMDP is upper bounded by the maximum entropy of its corresponding fully observable MDP. For a given controller $\pi$$\in$$\Pi(\mathbf{M})$ on a POMDP $\mathbf{M}$, we can construct, through \eqref{MDP_controller}, a controller $\pi'$$\in$$\Pi(\mathbf{M}_{fo})$ on the corresponding fully observable MDP $\mathbf{M}_{fo}$ which satisfies $Pr^{\pi}(s^{t+1} | s^{t})$$=$$Pr^{\pi'}(s^{t+1} | s^{t})$ for all $s^t$$\in$$\mathcal{SH}^t$ and $s^{t+1}$$\in$$\mathcal{SH}^{t+1}$. Then, for all $s^t$$\in$$\mathcal{SH}^t$, we have 
\begin{align*}
    \sup_{\pi\in\Pi(\mathbf{M})}H^{\pi}(S_{t+1}|S^{t}=s^{t})\leq \sup_{\pi\in\Pi(\mathbf{M}_{fo})}H^{\pi}(S_{t+1}|S^{t}=s^{t}).
\end{align*}

Informally, having access to the state history $s^t$, a controller $\pi'$$\in$$\Pi(\mathbf{M}_{fo})$ can attain an immediate reward $H^{\pi'}(S_{t+1}|S^{t}$$=$$s^{t})$ in \eqref{bellman_result} that is at least as high as the immediate reward achieved by a controller $\pi$$\in$$\Pi(\mathbf{M})$. Then, we have the following result as a consequence of Lemma \ref{writeValueFunction}.
{\setlength{\parindent}{0cm}
\noindent\begin{thm}\label{thm:POMDPbounded}
    For a POMDP $\mathbf{M}$ and a finite constant $N \in \mathbb{N}$,
    \begin{equation}\label{eq:POMDPbounded}
            \sup_{\pi \in \Pi(\mathbf{M})}H^{\pi}(S_1, S_2,\ldots, S_N)\leq \sup_{\pi \in \Pi(\mathbf{M}_{fo})}H^{\pi}(S_1, S_2,\ldots, S_N). 
    \end{equation}
\end{thm}}

The extension of Theorem \ref{thm:POMDPbounded} to infinite state sequences, i.e., to the case where $N$$\rightarrow$$\infty$, is rather straightforward. Since ${V}^{\pi}_{t,N}$, defined in \eqref{eq:defEntropyMDPsum}, is monotonically non-decreasing in $N$ for all $\pi$$\in$$\Pi(\mathbf{M})$, i.e., ${V}^{\pi}_{t,N+1}$$\geq$${V}^{\pi}_{t,N}$, we have 
\begin{align}\label{limsup}
    \sup_{\pi\in\Pi(\mathbf{M})}\lim_{N\rightarrow \infty} {V}_{t,N}^{\pi}(s^t)=\lim_{N\rightarrow \infty} \sup_{\pi\in\Pi(\mathbf{M})}{V}_{t,N}^{\pi}(s^t)
\end{align}
for all $s^t$$\in$$\mathcal{SH}^t$. Therefore, by taking the limits of both sides in \eqref{eq:POMDPbounded}, we conclude that the inequality between supremums still holds as $N$$\rightarrow$$\infty$. Finally, we conclude this section by noting that, with a slight modification of the statement of Lemma \ref{writeValueFunction}, it can be shown that all results presented in this section hold even if the future entropy rewards are discounted as in \eqref{objective2}.

\section{Entropy Maximization Over Finite-State Controllers}\label{section:Entropy Maximization Over Finite-State Controllers}

Optimal controllers solving the entropy maximization problems may, in general, use the complete system history to determine the next action to perform. To synthesize controllers in a POMDP, one can construct the so-called belief MDP whose (potentially infinite) states represent sufficient statistics of the system histories \cite{Kaelbling:1998:PAP:1643275.1643301,astrom}. However, such an approach is, in general, intractable as the number of states in the belief MDP grows exponentially with the length of the system histories. A common approach to overcome intractability is to restrict attention to finite-state controllers (FSCs) whose memory states represent (potentially insufficient) statistics of the system histories \cite{Poupart,meuleau1999solving, hansen}. Accordingly, in this section, we focus on FSCs with a fixed number of memory states and develop methods to synthesize locally optimal controllers within this restricted domain.
{\setlength{\parindent}{0cm}
\begin{definition}\label{def:FSC}
 For a POMDP $\mathbf{M}$, a $k$-\textit{finite-state controller} ($k$-FSC) is a tuple $\mathbf{C}$$=$$(\mathcal{Q},q_1,\gamma,\delta)$, where $\mathcal{Q}$$=$$\{q_1,q_2,\ldots,q_{k}\}$ is a finite set of memory states, $q_{1}$$\in$$\mathcal{Q}$ is the initial memory state, $\gamma$$:$$\mathcal{Q}$$\times$$\mathcal{Z}$$\rightarrow$$ \Delta(\mathcal{A})$ is a decision function and $\delta$$:$$\mathcal{Q}$$\times$$\mathcal{Z}$$\times$$\mathcal{A}$$\rightarrow$$\Delta(\mathcal{Q})$ is a memory transition function. We denote the collection of all $k$-FSCs by $\mathcal{F}_k(\mathbf{M})$.
\end{definition}}

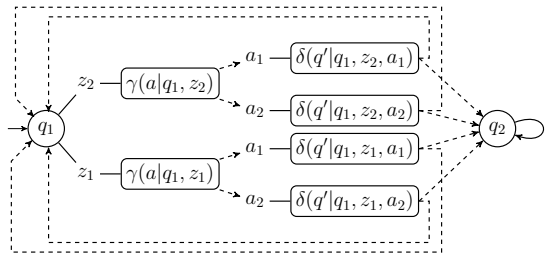
\begin{figure}[t!]\centering
\scalebox{0.55}{
\begin{tikzpicture}[->, >=stealth', auto, semithick, node distance=2cm]

     \tikzstyle{every state}=[fill=white,draw=black,thick,text=black,scale=1]

     \node[state,initial,initial text=] (q_1) {\Large{$q_1$}};
     \node        (z_1) [below right=5mm and 3mm of q_1] {\Large{$z_1$}};
     \node[draw=black, thick, rounded corners]   (gamma_1) [right=5mm of z_1] {\Large{$\gamma (a\lvert q_1,z_1)$}};
     \node   (a_11) [above right=-0.5 mm and 5mm of gamma_1] {\Large{$a_1$}};
     \node   (a_12) [below right=0 mm and 5mm of gamma_1] {\Large{$a_2$}};
     \node        (z_2) [above right=5mm and 3mm of q_1] {\Large{$z_2$}};
     \node[draw=black, thick, rounded corners]   (gamma_2) [right=5mm of z_2] {\Large{$\gamma(a\lvert q_1,z_2)$}};
     \node   (a_21) [above right=-0.5 mm and 5mm of gamma_2] {\Large{$a_1$}};
     \node   (a_22) [below right=0 mm and 5mm of gamma_2] {\Large{$a_2$}};
     \node[draw=black, thick, rounded corners]   (delta_11) [right= 5 mm of a_11] {\Large{$\delta(q'\lvert q_1,z_1,a_1)$}};  
     \node[draw=black, thick,  rounded corners]   (delta_12) [right= 5 mm of a_12] {\Large{$\delta(q'\lvert q_1,z_1,a_2)$}};  
     \node[draw=black, thick, rounded corners]   (delta_21) [right= 5 mm of a_21] {\Large{$\delta(q'\lvert q_1,z_2,a_1)$}};  
     \node[draw=black, thick,  rounded corners]   (delta_22) [right= 5 mm of a_22] {\Large{$\delta(q'\lvert q_1,z_2,a_2)$}};  
     \node[state] (q_2) [right=100mm of q_1]  {\Large{$q_2$}};

    \draw[thick,-] (q_1) -- (z_1);
    \draw[thick,-] (z_1) -- (gamma_1);
    \draw[thick,dashed, ->] (gamma_1) -- (a_11);
    \draw[thick,dashed, ->] (gamma_1) -- (a_12);
    \draw[thick,-] (a_11) -- (delta_11);
    \draw[thick,-] (a_12) -- (delta_12);
    \draw[thick,-] (q_1) -- (z_2);
    \draw[thick,-] (z_2) -- (gamma_2);
    \draw[thick,dashed, ->] (gamma_2) -- (a_21);
    \draw[thick,dashed, ->] (gamma_2) -- (a_22);
    \draw[thick,-] (a_21) -- (delta_21);
    \draw[thick,-] (a_22) -- (delta_22);
    \draw[thick, dashed, ->] (delta_11.0) -- (q_2.190);
    \draw[thick, dashed, ->] (delta_12.0) -- (q_2.220);
    \draw[thick, dashed, ->] (delta_21.0) -- (q_2.140);
    \draw[thick, dashed, ->] (delta_22.0) -- (q_2.170);
    \draw[thick,->] (q_2) to [out=380,in=340,looseness=8] (q_2);
    \draw[thick,dashed,->] (delta_12) -| ++(1.8,-1) -| ++(-9.2,2.3);
    \draw[thick,dashed,->] (delta_11) -| ++(2.1,-2.5) -| ++(-10.4,2.2) -- ++(0.5,0.5);
    \draw[thick,dashed,->] (delta_21) -| ++(1.8,1) -| ++(-9.2,-2.3);
    \draw[thick,dashed,->] (delta_22) -| ++(2.1,2.5) -| ++(-10.3,-2.2) -- ++(0.45,-0.45);
\end{tikzpicture}}
\caption{An illustration of finite-state controllers. In a memory state $q$, the agent receives an observation $z$, chooses an action $a$ based on the decision function $\gamma(a\lvert q,z)$, and transitions to a memory state $q'$ based on the transition function $\delta(q'\lvert q,z,a)$. The functions $\gamma$ and $\delta$ are design variables, and their outcomes are indicated with dashed lines. }\label{FSC_illustration}
\end{figure}

In Fig. \ref{FSC_illustration}, we present an illustration of $k$-FSCs. As can be seen in the figure, in a memory state representing certain statistics of system histories in a POMDP, the agent receives an observation, makes a decision based on the function $\gamma$, and updates its memory state based on the function $\delta$.

For a POMDP, a $k$-FSC induces a Markov chain (MC) which is an MDP with a single available action, i.e., $\lvert \mathcal{A} \rvert$$=$$1$. It is shown in \cite{junges2017permissive} that the set of all MCs that can be induced by a $k$-FSC is the set of all well-defined instantiations of a specific parametric MC (pMC). Therefore, without loss of generality, one can work on that specific pMC to synthesize an instantiation which corresponds to the MC induced by an entropy-maximizing FSC. In the following sections, we first provide formal definitions of pMCs and their instantiations. We then reformulate the entropy maximization problem over $k$-FSCs as another optimization problem over pMCs.
{\setlength{\parindent}{0cm}
\begin{remark}
Any given instance of the finite horizon entropy maximization problem can be reduced to an instance of the infinite horizon entropy maximization problem in time polynomial in the size of the POMDP and the decision horizon $N$. For brevity, we establish our results only for the infinite horizon entropy maximization problem and provide the details of the aforementioned reduction in Appendix \ref{section:Appendix2}. 
\end{remark}}


\subsection{Parametric Markov Chains}
We develop solutions to entropy maximization problems through the use of parametric Markov chains.
 {\setlength{\parindent}{0cm}\noindent
\begin{definition}\label{def:pMC}
    For a POMDP $\mathbf{M}$ and a constant $k$$\in$$\mathbb{N}$, the \textit{induced parametric Markov chain} (pMC) is a tuple $\mathbf{D}_{\mathbf{M},k}$$=$$ (\mathcal{S}_{\mathbf{M},k},s_{I,\mathbf{M},k},\mathcal{V}_{\mathbf{M},k},P_{\mathbf{M},k}, {R}_{\mathbf{M},k})$ where entries are as follows;
    \begin{itemize}
        \item $\mathcal{S}_{\mathbf{M},k} = \mathcal{S} \times \{1,2,...,k\}$ is the finite set of states.
        \item $s_{I,\mathbf{M},k} = \langle s_{I}, 1 \rangle$ is the initial state.
        \item $\!\begin{aligned}[t]
                    \mathcal{V}_{\mathbf{M},k} = \{\gamma_{a}^{q,z}| & z \in \mathcal{Z},q \in \mathcal{Q}, a \in \mathcal{A} \} \nonumber \\
                    & \cup \{\delta_{q'}^{q,z,a}|z \in \mathcal{Z}, q,q' \in \mathcal{Q}, a \in \mathcal{A} \} \nonumber
                \end{aligned}$\\
                is the finite set of parameters.
        \item $P_{\mathbf{M},k}$ $:$ $\mathcal{S}_{\mathbf{M},k}$$\rightarrow$$\Delta(\mathcal{S}_{\mathbf{M},k}) $ 
        is a transition function such that $P_{\mathbf{M},k}(\langle s',q' \rangle\ | \ \langle s,q \rangle):= \sum_{a\in \mathcal{A}}\overline{P}(\langle s',q' \rangle \ | \ \langle s,q \rangle,a)$ for all $\langle s,q \rangle,\langle s',q'\rangle$$\in$$\mathcal{S}_{\mathbf{M},k}$ where $\overline{P}$ $:$ $\mathcal{S}_{\mathbf{M},k}$$\times$$\mathcal{A}$$\rightarrow$$ \Delta(\mathcal{S}_{\mathbf{M},k})$ is a mapping such that
         \begin{flalign}\label{trans_pmc}
        \overline{P}(\langle s',q' \rangle \ | \ \langle s,q \rangle,a) := \sum_{z \in \mathcal{Z}} {O}_{s,z} \, {P}_{s,a,s'} \, \gamma_{a}^{q,z} \, \delta_{q'}^{q,z,a}.&& \raisetag{22pt}
    \end{flalign}
        \item ${R}_{\mathbf{M},k}(\langle s,q \rangle,a):={R}(s,a)$ for all $s$$\in$$\mathcal{S}$, $q\in \mathcal{Q}$, and $a$$\in$$\mathcal{A}$.
 \end{itemize}
\end{definition}}


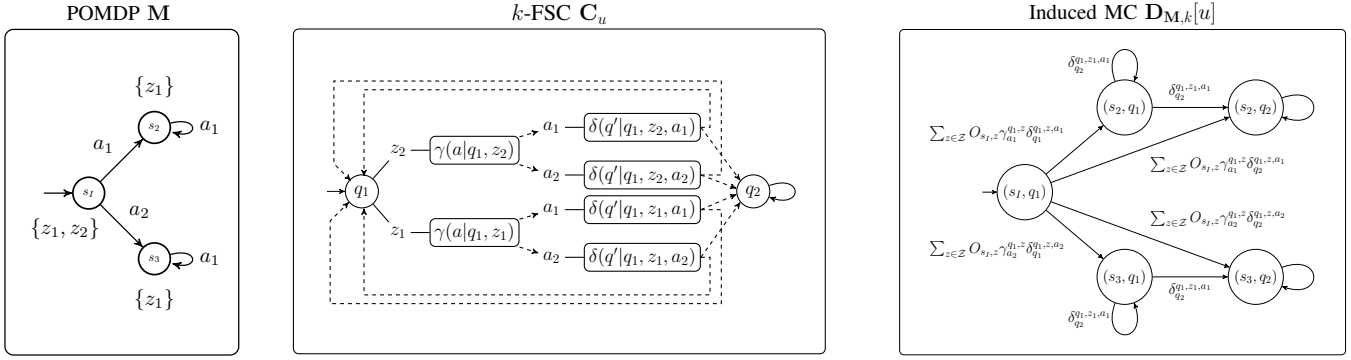
\begin{figure*}[t!]
\begin{subfigure}{.2\linewidth}
\scalebox{0.8}{
\begin{tikzpicture}[->, >=stealth', auto, semithick, node distance=2cm]

    \tikzstyle{every state}=[fill=white,draw=black,thick,text=black,scale=0.6]

    \node[state,initial,initial text=] (s_1) {$s_I$};
    \node[] (summation) [right= 22.5 mm of s_1] {};
    \node[] (summation) [right= 120.5 mm of s_1] {};
    \node[state] (s_2) [above right =10mm of s_1]  {$s_2$};
    \node[] (s_1z) [below left =1mm and -5mm of s_1]  {$\{z_1,z_2\}$};
    \node[] (s_2z) [above =1mm of s_2]  {$\{z_1\}$};
    \node[state] (s_3) [below right =10mm of s_1]  {$s_3$};
    \node[] (s_3z) [below =1mm of s_3]  {$\{z_1\}$};

      \draw[rounded corners] ($(s_1.north west)+(-1.2,2.53)$)  rectangle    ($(s_1.south     east)+(2.3,-2.5)$); 
    \node[anchor=south, above left = 0.1 mm and 41 mm] at (current bounding box.north) {POMDP $\mathbf{M}$};

\path
(s_1)	 edge     node{$a_1$}     (s_2)
(s_1)	 edge     node{$a_2$}     (s_3)
(s_2)  edge  [loop right=10]    node{$a_1$}     (s_2)
(s_3)  edge  [loop right=10]    node{$a_1$}     (s_3);
\end{tikzpicture}}
\end{subfigure}~ 
\begin{subfigure}{.2\linewidth}
\scalebox{0.5}{
\begin{tikzpicture}[->, >=stealth', auto, semithick, node distance=2cm]

     \tikzstyle{every state}=[fill=white,draw=black,thick,text=black,scale=1]

     \node[state,initial,initial text=] (q_1) {\Large{$q_1$}};
     \node        (z_1) [below right=5mm and 3mm of q_1] {\Large{$z_1$}};
     \node[draw=black, thick, rounded corners]   (gamma_1) [right=5mm of z_1] {\Large{$\gamma(a \lvert q_1,z_1)$}};
     \node   (a_11) [above right=-0.5 mm and 5mm of gamma_1] {\Large{$a_1$}};
     \node   (a_12) [below right=0 mm and 5mm of gamma_1] {\Large{$a_2$}};
     \node        (z_2) [above right=5mm and 3mm of q_1] {\Large{$z_2$}};
     \node[draw=black, thick, rounded corners]   (gamma_2) [right=5mm of z_2] {\Large{$\gamma(a \lvert q_1,z_2)$}};
     \node   (a_21) [above right=-0.5 mm and 5mm of gamma_2] {\Large{$a_1$}};
     \node   (a_22) [below right=0 mm and 5mm of gamma_2] {\Large{$a_2$}};
     \node[draw=black, thick, rounded corners]   (delta_11) [right= 5 mm of a_11] {\Large{$\delta(q' \lvert q_1,z_1,a_1)$}};  
     \node[draw=black, thick, rounded corners]   (delta_12) [right= 5 mm of a_12] {\Large{$\delta(q'\lvert q_1,z_1,a_2)$}};  
     \node[draw=black, thick,rounded corners]   (delta_21) [right= 5 mm of a_21] {\Large{$\delta(q'\lvert q_1,z_2,a_1)$}};  
     \node[draw=black, thick, rounded corners]   (delta_22) [right= 5 mm of a_22] {\Large{$\delta(q'\lvert q_1,z_2,a_2)$}};  
     \node[state] (q_2) [right=95mm of q_1]  {\Large{$q_2$}};
     
      \draw[thick,rounded corners] ($(q_1.north west)+(-1.5,4)$)  rectangle    ($(q_1.south     east)+(12,-4)$); 
    \node[anchor=south] at (current bounding box.north) {\LARGE $k$-FSC $\mathbf{C}_u$};

    \draw[thick,-] (q_1) -- (z_1);
    \draw[thick,-] (z_1) -- (gamma_1);
    \draw[thick,dashed, ->] (gamma_1) -- (a_11);
    \draw[thick,dashed, ->] (gamma_1) -- (a_12);
    \draw[thick,-] (a_11) -- (delta_11);
    \draw[thick,-] (a_12) -- (delta_12);
    \draw[thick,-] (q_1) -- (z_2);
    \draw[thick,-] (z_2) -- (gamma_2);
    \draw[thick,dashed, ->] (gamma_2) -- (a_21);
    \draw[thick,dashed, ->] (gamma_2) -- (a_22);
    \draw[thick,-] (a_21) -- (delta_21);
    \draw[thick,-] (a_22) -- (delta_22);
    \draw[thick, dashed, ->] (delta_11.0) -- (q_2.190);
    \draw[thick, dashed, ->] (delta_12.0) -- (q_2.220);
    \draw[thick, dashed, ->] (delta_21.0) -- (q_2.140);
    \draw[thick, dashed, ->] (delta_22.0) -- (q_2.170);
    \draw[thick,->] (q_2) to [out=380,in=340,looseness=8] (q_2);
    \draw[thick,dashed,->] (delta_12) -| ++(1.8,-1) -| ++(-9.2,2.3);
    \draw[thick,dashed,->] (delta_11) -| ++(2.1,-2.5) -| ++(-10.4,2.2) -- ++(0.5,0.5);
    \draw[thick,dashed,->] (delta_21) -| ++(1.8,1) -| ++(-9.2,-2.3);
    \draw[thick,dashed,->] (delta_22) -| ++(2.1,2.5) -| ++(-10.3,-2.2) -- ++(0.45,-0.45);
\end{tikzpicture}}
\end{subfigure}~\hspace{42mm}
\begin{subfigure}{.2\linewidth}
\scalebox{0.4}{
\begin{tikzpicture}[->, >=stealth', auto, semithick, node distance=2cm]

     \tikzstyle{every state}=[fill=white,draw=black,thick,text=black,scale=1]

     \node[state,initial,initial text=] (s_1_q_1) {\Large{$(s_I,q_1)$}};
     \node[state] (s_2_q_1) [above right =15mm and 20 mm of s_1_q_1]  {\Large{$(s_2,q_1)$}};
     \node[state] (s_3_q_1) [below right =15mm and 20 mm of s_1_q_1]  {\Large{$(s_3,q_1)$}};
     \node[state] (s_2_q_2) [right =25mm of s_2_q_1]  {\Large{$(s_2,q_2)$}};
     \node[state] (s_3_q_2) [right =25mm of s_3_q_1]  {\Large{$(s_3,q_2)$}};
     \draw[thick,->] (s_1_q_1) -- node[above left = 1 mm and 2 mm]{\Large $\sum_{z\in \mathcal{Z}}O_{s_I,z}\gamma^{q_1,z}_{a_1}\delta^{q_1,z,a_1}_{q_1}$} (s_2_q_1);
     \draw[thick,->] (s_1_q_1) -- node[below left = 1 mm and 2 mm]{\Large $\sum_{z\in \mathcal{Z}}O_{s_I,z}\gamma^{q_1,z}_{a_2}\delta^{q_1,z,a_2}_{q_1}$} (s_3_q_1);
     \draw[thick,->] (s_1_q_1) -- node[above right = 1 mm and 2 mm]{\Large $\sum_{z\in \mathcal{Z}}O_{s_I,z}\gamma^{q_1,z}_{a_2}\delta^{q_1,z,a_2}_{q_2}$} (s_3_q_2);
     \draw[thick,->] (s_1_q_1) -- node[below right = 1 mm and 2 mm]{\Large $\sum_{z\in \mathcal{Z}}O_{s_I,z}\gamma^{q_1,z}_{a_1}\delta^{q_1,z,a_1}_{q_2}$} (s_2_q_2);
     \draw[thick,->] (s_2_q_1) -- node[above = 1 mm]{\Large $\delta^{q_1,z_1,a_1}_{q_2}$} (s_2_q_2);
     \draw[thick,->] (s_3_q_1) -- node[below = 1 mm ]{\Large $\delta^{q_1,z_1,a_1}_{q_2}$} (s_3_q_2);

    \draw[thick,->] (s_2_q_1) to [out=110,in=70,looseness=6] node[below left = 1 mm and 3.2 mm]{\Large $\delta^{q_1,z_1,a_1}_{q_2}$} (s_2_q_1);
    \draw[thick,->] (s_2_q_2) to [out=380,in=340,looseness=6] (s_2_q_2);
    \draw[thick,->] (s_3_q_1) to [out=250,in=290,looseness=6] node[above left = 1 mm and 3.2 mm]{\Large $\delta^{q_1,z_1,a_1}_{q_2}$} (s_3_q_1);
    \draw[thick,->] (s_3_q_2) to [out=380,in=340,looseness=6] (s_3_q_2);

      \draw[thick,rounded corners] ($(s_1_q_1.north west)+(-3.5,4.75)$)  rectangle    ($(s_1_q_1.south     east)+(10,-4.75)$);
    \node[anchor=south] at (current bounding box.north) {\huge Induced MC $\mathbf{D}_{\mathbf{M},k}[u]$};
 
\end{tikzpicture}}
\end{subfigure}
\caption{ An illustration of an MC $\mathbf{D}_{\mathbf{M},k}[u]$ (right) induced from a POMDP $\mathbf{M}$ (left) by a $k$-FSC $\mathbf{C}_u$ (middle). 
Note that, the number of states in $\mathbf{D}_{\mathbf{M},k}[u]$ is larger than $\mathbf{M}$ due to the stochasticity in the memory transition function $\delta$. }\label{pMC_illustration}
\end{figure*}


An MC can be obtained from an induced pMC by instantiating the parameters $\mathcal{V}_{\mathbf{M},k}$ in a way that the resulting transition function is well-defined. Formally, a \textit{well-defined instantiation} for $\mathcal{V}_{\mathbf{M},k}$ is a function $u$ $:$ $ \mathcal{V}_{\mathbf{M},k} $$\rightarrow$$[0,1]$ such that, for all $a\in\mathcal{A}$, $q\in \mathcal{Q}$, and $z\in\mathcal{Z}$,
\begin{align*}
    \sum_{a\in\mathcal{A}}u(\gamma_a^{q,z})=1\  \ \ \text{and} \ \ \ \sum_{q'\in \mathcal{Q}}u(\delta_{q'}^{q,z,a})=1.
\end{align*}

Applying a well-defined instantiation $u$ to the induced pMC $\mathbf{D}_{\mathbf{M},k}$, denoted $\mathbf{D}_{\mathbf{M},k}[u]$, replaces each parametric transition probability $P_{\mathbf{M},k}$ by $P^u_{\mathbf{M},k}$. It is straightforward to verify that $\mathbf{D}_{\mathbf{M},k}[u]$ is an MC.
Let $\Upsilon_{\mathbf{M},k}$ denote the set of all well-defined instantiations for a pMC $\mathbf{D}_{\mathbf{M},k}$. For an induced pMC $\mathbf{D}_{\mathbf{M},k}$, every instantiation $u$$\in$$\Upsilon_{\mathbf{M},k}$ describes a k-FSC $\mathbf{C}_u$$\in$$\mathcal{F}_k(\mathbf{M})$ \cite{junges2017permissive}. Thus, we can synthesize all admissible MCs that can be induced from a POMDP $\mathbf{M}$ by a k-FSC $\mathbf{C}$$\in$$\mathcal{F}_k(\mathbf{M})$ through well-defined instantiations $u$ over $\mathcal{V}_{\mathbf{M},k}$. In Fig. \ref{pMC_illustration}, we provide a simple example to illustrate the derivation of the instantiation $\mathbf{D}_{\mathbf{M},k}[u]$ from a given POMDP $\mathbf{M}$ and a $k$-FSC $\mathbf{C}_u$.


\subsection{A Reformulation of Entropy Maximization Problem over Parametric Markov Chains }
Recall that we are interested in maximizing the entropy of the state sequence of a given POMDP $\mathbf{M}$. As can be seen from Fig. \ref{pMC_illustration}, the number of states that are reachable from the initial state of a pMC $\mathbf{D}_{\mathbf{M},k}$ is, in general, larger than that of the POMDP $\mathbf{M}$. It is known \cite{Cover} that the maximum entropy of a random variable increases as the cardinality of its support increases. Hence, it is possible to have a well-defined instantiation $\mathbf{D}_{\mathbf{M},k}[u]$ whose entropy of state sequences is higher than the maximum entropy of the state sequences of the POMDP $\mathbf{M}$. This observation implies that, in general, the maximum entropy of a POMDP $\mathbf{M}$ is \textit{not} an upper bound on the maximum entropy of the induced pMC $\mathbf{D}_{\mathbf{M},k}$.

The maximum entropy of $\mathbf{D}_{\mathbf{M},k}$ is, in general, higher than that of $\mathbf{M}$ due to the stochasticity introduced to the process by the parameters $\delta^{q,z,a}_{q'}$. To synthesize an entropy-maximizing $k$-FSC for a POMDP $\mathbf{M}$ using the induced pMC $\mathbf{D}_{\mathbf{M},k}$, we impose restrictions on the memory transition function of the $k$-FSCs. For each memory state $q$$\in$$\mathcal{Q}$ in a given $k$-FSC, let 
\begin{align*}
    Succ(q):=\{q'\in\mathcal{Q} : \delta(q' | q,z,a)>0, z\in \mathcal{Z}, a \in \mathcal{A}\}.
\end{align*}
\vspace{-0.5cm}
{\setlength{\parindent}{0cm}
\begin{definition}\label{deterministicFSC}
A \textit{deterministic $k$-FSC} $\mathbf{C}$$=$$(\mathcal{Q},q_1,\gamma,\delta)$ is a $k$-FSC such that for all $q$$\in$$\mathcal{Q}$, $|Succ(q)| = 1$. We denote the set of all deterministic k-FSCs by $\mathcal{F}_{k}^{det}(\mathbf{M})$.
\end{definition}}

For a $k$-FSC $\mathbf{C}$$\in$$\mathcal{F}_{k}^{det}(\mathbf{M})$, let $u_{\mathbf{C}}$$:$$V_{\mathbf{M},k}$$\rightarrow$$\mathbb{R}$ be the corresponding instantiation of the induced pMC $\mathbf{D}_{\mathbf{M},k}$ such that
\begin{align*}
    u_{\mathbf{C}}(\gamma_{a}^{q,z}):=\gamma(a|q,z) \ \  \text{and} \ \  u_{\mathbf{C}}(\delta_{q'}^{q,z,a}):=\delta(q'|q,z,a).
\end{align*}
 Moreover, let $\Upsilon_{\mathbf{M},k}^{det}$ denote the set of all well-defined instantiations $u_{\mathbf{C}}$ that corresponds to a deterministic $k$-FSC $\mathbf{C}$. Noting that $\mathbf{D}_{\mathbf{M},k}[u_{\mathbf{C}}]$ is a stochastic process, we denote its sequence of states by $({S}_{\mathbf{M},k,1}, S_{\mathbf{M},k,2},\ldots)$. Moreover, for a given state ${S}_{\mathbf{M},k,t-1}$, we denote the one-step entropy of $\mathbf{D}_{\mathbf{M},k}[u_{\mathbf{C}}]$ by $H^{u_{\mathbf{C}}}(S_{\mathbf{M},k,t} | S_{\mathbf{M},k,t-1})$. 
\begin{prop}\label{pmc_prop}
For a given POMDP $\mathbf{M}$, a controller $\mathbf{C}$$\in$$\mathcal{F}_{k}^{det}(\mathbf{M})$, and constants $t,k$$\in$$\mathbb{N}$, we have
\begin{align}\label{state_hist_to_pMC}
H^{\mathbf{C}}(S_{t} | S^{t-1})=  H^{u_{\mathbf{C}}}(S_{\mathbf{M},k,t} | S_{\mathbf{M},k, t-1}).
\end{align}
\end{prop}

Proposition \ref{pmc_prop} shows that the local (one-step) entropy gained in a POMDP $\mathbf{M}$ under a deterministic $k$-FSC $\mathbf{C}$ is equal to the local entropy gained in $\mathbf{D}_{\mathbf{M},k}[u_{\mathbf{C}}]$. Note that, since the memory states $q$ are explicitly represented in the states $\langle s,q\rangle$ of $\mathbf{D}_{\mathbf{M},k}[u_{\mathbf{C}}]$, one-step entropy in $\mathbf{D}_{\mathbf{M},k}[u_{\mathbf{C}}]$ depends only on the state occupied in the previous step.

Proposition \ref{pmc_prop}, together with the definition of the induced pMC, implies that, for any $\mathbf{C}$$\in$$\mathcal{F}_{k}^{det}(\mathbf{M})$,

\begin{subequations}
\begin{flalign}\label{pp11}
    \mathbf{C}\in\arg&\max_{\mathbf{C}'\in\mathcal{F}_{k}^{det}(\mathbf{M})}\ \ \sum_{t=2}^{\infty}\beta^{t-2}H^{\mathbf{C}'}(S_{t} | S^{t-1})&&\\ 
    &\text{subject to:} \ \ \label{pp12} \mathbb{E}^{\mathbf{C}'}\Big[\sum_{t=1}^{\infty}\beta^{t-1}{R}(S_t,A_t)\Big]\geq \Gamma &&\\ \nonumber \\ 
   & \qquad\qquad \qquad \text{if and only if}\nonumber
\end{flalign}
\end{subequations}
\begin{subequations}
\begin{flalign}\label{pp13}
   u_{\mathbf{C}}\in\arg& \max_{u\in\Upsilon_{\mathbf{M},k}^{det}}\ \ \sum_{t=2}^{\infty}\beta^{t-2}H^{u}(S_{\mathbf{M},k,t} | S_{\mathbf{M},k,t-1})&&\\  \label{pp14}
    & \text{subject to:} \ \ \mathbb{E}^{u}\Big[\sum_{t=1}^{\infty}\beta^{t-1}{R}(S_{\mathbf{M},k,t},A_t)\Big]\geq \Gamma. \raisetag{22pt}&&
\end{flalign}
\end{subequations}

The following result is due to the fact that $\mathcal{F}_k^{det}(\mathbf{M})$$\subseteq$$\Pi(\mathbf{M})$ and shows that the maximum entropy of the pMC induced from FSCs with deterministic memory transitions is upper bounded by that of the corresponding POMDP.

{\setlength{\parindent}{0cm}
\begin{Corollary}
Let $G_1^{\star}$ and $G_2^{\star}$ be the optimal values of the problems given in \eqref{objective2}-\eqref{constraint2} and \eqref{pp13}-\eqref{pp14}, respectively. Then, we have $G_1^{\star}$$\geq$$G_2^{\star}$.
\end{Corollary}}
{\setlength{\parindent}{0cm}
\begin{remark}
We emphasize that the equivalence between the problems \eqref{pp11}-\eqref{pp12} and \eqref{pp13}-\eqref{pp14} does not hold if we consider the entropy maximization problem over the set $\mathcal{F}_k(\mathbf{M})$ instead of the set  $\mathcal{F}^{det}_k(\mathbf{M})$. Note that for a POMDP $\mathbf{M}$ under a controller $\mathbf{C}$$\in$$\mathcal{F}_k(\mathbf{M})$, we have $\lvert Succ(s)\rvert$$\leq$$\lvert \mathcal{S}\rvert$ for all $s$$\in$$\mathcal{S}$, whereas on the corresponding induced pMC, it is possible to have $\lvert Succ(\langle s,q\rangle)\lvert$$>$$\lvert \mathcal{S} \rvert$ as can be seen in Fig. \ref{pMC_illustration}. The entropy of a distribution with a support of size $\lvert Succ(\langle s,q\rangle)\lvert$ can be made larger than the entropy of a distribution with a support of size $\lvert Succ(s)\rvert$. Hence, if one performs the maximization over the set $\mathcal{F}_k(\mathbf{M})$, the maximum entropy of the induced pMC $\mathbf{D}_{\mathbf{M},k}$ is not necessarily equal to the maximum entropy of the corresponding POMDP $\mathbf{M}$.
\end{remark}}

\subsection{Finite-State Controller Synthesis: Optimization Problem}\label{FSCsynthesis}

We now present a nonlinear optimization problem to synthesize a deterministic $k$-FSC that maximizes the entropy of a POMDP over all deterministic $k$-FSCs.

Recall that for a POMDP $\mathbf{M}$ and a constant $k$$>$$0$, the induced pMC $\mathbf{D}_{\mathbf{M},k}$ represents all possible MCs that can be induced from $\mathbf{M}$ by a $k$-FSC. Moreover, Proposition \ref{pmc_prop} implies that the maximum entropy of $\mathbf{D}_{\mathbf{M},k}$ is equal to the maximum entropy of $\mathbf{M}$ \textit{if one restricts attention to k-FSCs with deterministic memory transitions}. In what follows, we formulate an optimization problem to synthesize a well-defined instantiation $u$ for the pMC $\mathbf{D}_{\mathbf{M},k}$ such that the entropy of the MC $\mathbf{D}_{\mathbf{M},k}[u]$ is maximized over all MCs $\mathbf{D}_{\mathbf{M},k}[u']$ for which $P^{u'}_{\mathbf{M},k}(\langle s',q' \rangle\ | \ \langle s,q \rangle)$$>$$0$ for a single $q'$$\in$$\mathcal{Q}$. 

To restrict the search space to FSCs with deterministic memory transitions, we introduce the following constraints 
\begin{align*}
    u(\delta_{q'}^{q,z,a})\in\{0,1\} \ \ \text{and}\ \  \sum_{z\in\mathcal{Z}}\sum_{a\in\mathcal{A}}u(\delta_{q'}^{q,z,a})\in\{0,\lvert \mathcal{Z}\rvert\lvert \mathcal{A}\rvert\}. 
\end{align*}
Intuitively, the above constraints ensure the  transition to a single successor memory state regardless of the received observation and the taken action. We note that, for a given memory state pair $(q,q')$, the second integer constraint can be implemented as $\lvert\mathcal{Z}\rvert\lvert \mathcal{A}\lvert$ equality constraints. Finally, the above constraints do not prevent the agent from randomizing its actions. The agent can still randomize its actions at a given state $s$$\in$$\mathcal{S}$ by instantiating the parameters $\gamma_{a}^{q,z}$ appropriately.

For notational simplicity, let ${\bf{s}}$ denote an arbitrary state $\langle s,t\rangle$$\in$$\mathcal{S}_{\mathbf{M},k}$. Let $L^{u}$$:$$\mathcal{S}_{\mathbf{M},k}$$\rightarrow$$\mathbb{R}$ be the \textit{local entropy} function such that, for all ${\bf{s}}$$\in$$\mathcal{S}_{\mathbf{M},k}$,
\begin{align}
L^u({\bf{s}}):=-\sum_{{\bf{s}'}\in \mathcal{S}_{\mathbf{M},k}}P^u_{\mathbf{M},k}({\bf{s'}}|{\bf{s}})\log P^u_{\mathbf{M},k}({\bf{s'}}|{\bf{s}}).
\end{align}
Note that we have $L^u({\bf{s}})$$=$$H^{u}(S_{\mathbf{M},k,t} | S_{\mathbf{M},k, t-1}$$=$${\bf{s}})$ for any $t$$\in$$\mathbb{N}$. Hence, $L^u({\bf{s}})$ corresponds to the one-step entropy reward gained in the MC $\mathbf{D}_{\mathbf{M},k}[u]$ from the state $\bf{s}$. Recalling the equivalence given in \eqref{state_hist_to_pMC}, the local entropy function $L^u$ allows us to transfer the results of Section \ref{section:An Upper Bound on Maximum Entropy}, which are obtained for a POMDP $\mathbf{M}$, to the pMC $\mathbf{D}_{\mathbf{M},k}$. Specifically, by slightly modifying the statement of Lemma \ref{writeValueFunction} and defining variables $\nu$$\in$$\mathbb{R}^{\lvert \mathcal{S}_{\mathbf{M},k}\rvert}$, it can be shown that the maximum entropy \eqref{pp13} of $\mathbf{D}_{\mathbf{M},k}$ is the unique fixed-point of the system of equations
\begin{flalign}
  &  \nu({\bf{s}}) = \max_{u\in \Upsilon_{\mathbf{M},k}}\Bigg\{L^{u}({\bf{s}}) + \beta \sum_{{\bf{s'}} \in \mathcal{S}_{\mathbf{M},k}}P_{\mathbf{M},k}^{u} ({\bf{s'}}| {\bf{s}})\nu({\bf{s'}}) \Bigg\}
    \end{flalign}
and equal to $\nu({\bf{s}}_I)$$:=$$\nu(s_{I,\mathbf{M},k})$. Hence, the maximum entropy \eqref{pp13} of $\mathbf{D}_{\mathbf{M},k}$ can be computed by finding the maximum $\nu({\bf{s}}_I)$ that satisfies 
\begin{flalign*}
  &  \nu({\bf{s}}) \leq L^{u}({\bf{s}}) + \beta \sum_{{\bf{s}'} \in \mathcal{S}_{\mathbf{M},k}}P_{\mathbf{M},k}^{u} ({\bf{s}'}| {\bf{s}})\nu({\bf{s}'})  \ \forall {\bf{s}}\in \mathcal{S}_{\mathbf{M},k}.
    \end{flalign*}
    In the above inequality, both $P_{\mathbf{M},k}^{u} ({\bf{s}'}| {\bf{s}})$ and $\nu({\bf{s}'})$ are variables. Hence, standard methods, e.g., value iteration, cannot be used to compute $\nu({\bf{s}}_I)$; instead, one needs to solve a nonlinear optimization problem, which we present shortly, for the computation of $\nu({\bf{s}}_I)$. 
    
    Let ${R}^u$$:$$\mathcal{S}_{\mathbf{M},k}$$\rightarrow$$\mathbb{R}$ be the expected immediate rewards on $\mathbf{D}_{\mathbf{M},k}$ such that, for all ${\bf{s}}$$\in$$\mathcal{S}_{\mathbf{M},k}$,
\begin{align}
    {R}^u({\bf{s}}):=\sum_{{\bf{s}'} \in \mathcal{S}_{\mathbf{M},k}}\sum_{a\in\mathcal{A}}\overline{P}^{u} ({\bf{s}'}| {\bf{s}},a){R}({\bf{s}'},a)
\end{align}
where $\overline{P}^{u}$$:$$S_{\mathbf{M},k}$$\times$$ \mathcal{A}$$\rightarrow$$\Delta(\mathcal{S}_{\mathbf{M},k})$ is defined by replacing parameters $\gamma_a^{q,z}$ and $\delta^{q,z,a}_{q'}$ in \eqref{trans_pmc} with their corresponding instantiations $u(\gamma_a^{q,z})$ and $u(\delta^{q,z,a}_{q'})$. Then, the problem in \eqref{pp13}-\eqref{pp14} can be formulated as a nonlinear optimization problem (NLP) as follows:
\begin{subequations}
    \begin{align}
       & \underset{\nu,u, \eta}{\text{maximize}} \qquad  \nu({\bf{s}}_{I})&& \label{objectiveFunction} \\
    &    \text{subject to:}  &&\nonumber \\ \label{first_trilinear}
       &   \nu({\bf{s}}) \leq L^{u}({\bf{s}}) + \beta \sum_{{\bf{s}'} \in \mathcal{S}_{\mathbf{M},k}}P_{\mathbf{M},k}^{u} ({\bf{s}'}| {\bf{s}})\nu({\bf{s}'})  \ \forall {\bf{s}}\in \mathcal{S}_{\mathbf{M},k}\raisetag{26pt}  &&\\ 
      &  \eta({\bf{s}}) \leq {R}^u({\bf{s}})+\beta\sum_{{\bf{s}}' \in \mathcal{S}_{\mathbf{M},k}}P_{\mathbf{M},k}^{u} ({\bf{s}}'| {\bf{s}})\eta({\bf{s}}') \ \forall {\bf{s}} \in \mathcal{S}_{\mathbf{M},k}  \label{reachabilityProb}  \raisetag{26pt} &&\\
        &  \eta({\bf{s}}_{I}) \geq \Gamma  \label{initialStateProb}&&\\
&u(\delta_{q'}^{q,z,a})\in\{0,1\}, \ \ \sum_{q'\in Q} u(\delta_{q'}^{q,z,a})=1 \label{integer_cons}\\
& 0 \leq u(\gamma_{a}^{q,z})\leq 1,\ \  \ \ \sum_{a\in\mathcal{A}} u(\gamma_{a}^{q,z})=1\\
&\sum_{z\in\mathcal{Z}}\sum_{a\in\mathcal{A}}u(\delta_{q'}^{q,z,a})\in\{0,\lvert \mathcal{Z}\rvert\lvert \mathcal{A}\rvert\}\label{last_cons}.
    \end{align}
\end{subequations}

In the above optimization problem, the variable $\eta({\bf{s}})$ denotes the expected reward collected by starting from the state ${\bf{s}}$$\in$$\mathcal{S}_{\mathbf{M},k}$. It follows from \cite{puterman2014markov} that the constraint \eqref{initialStateProb} ensures that a solution $u^{\star}$ to the above problem collects an expected total reward exceeding the threshold $\Gamma$.

Recall that the transition function $P_{\mathbf{M},k}^{u}$ of the MC $\mathbf{D}_{\mathbf{M},k}[u]$, which results from the instantiation $u$ of the pMC $\mathbf{D}_{\mathbf{M},k}$, is given by $P_{\mathbf{M},k}^{u} ({\bf{s}'}| {\bf{s}})$$=$$\sum_{a\in \mathcal{A}} \overline{P}^u({\bf{s}'}| {\bf{s}},a)$ where 
\begin{align}\label{trilinear_expanded}
    \overline{P}^u({\bf{s}'}| {\bf{s}},a) := \sum_{z \in \mathcal{Z}} {O}_{s,z} \, {P}_{s,a,s'} \, u(\gamma_{a}^{q,z}) \, u(\delta_{q'}^{q,z,a}),
\end{align}
${\bf{s}}$$=$$\langle s,q \rangle$, and ${\bf{s}'}$$=$$\langle s',q' \rangle$. Therefore, the problem in \eqref{objectiveFunction}-\eqref{last_cons} involves nonlinear constraints in \eqref{first_trilinear} where three variables are multiplied with each other. Even though certain relaxation techniques, \textit{e.g.}, McCormick envelopes \cite{McCormick1976}, can be used to replace the constraints in \eqref{first_trilinear} with specific bilinear constraints, finding an optimal solution to the resulting NLP remains as a challenge due to binary constraints in \eqref{integer_cons}. One can employ branch-and-bound algorithms \cite{bertsekas1997nonlinear} to obtain a solution to the problem in \eqref{objectiveFunction}-\eqref{last_cons}, but unfortunately, such algorithms scale poorly with the size of the problem instances. 

For practical purposes, instead of computing a globally optimal solution, one can aim to obtain a locally optimal solution to the problem in \eqref{objectiveFunction}-\eqref{last_cons} \textit{after setting the instantiation $u(\delta_{q'}^{q,z,a})$ of memory transitions to a constant}. In the next section, we provide a technique to obtain a local optimal solution to the problem in \eqref{objective2}-\eqref{constraint2} over all $k$-FSCs with a specific deterministic memory transition function.

\subsection{Finite-State Controller Synthesis: A Solution Approach}\label{FSCsynthesis_solution}
In this section, we consider the entropy maximization problems over $k$-FSCs with a specific deterministic transition function and present an algorithm to synthesize a controller which locally maximizes the entropy of a given POMDP.

We first set the variables $u(\delta_{q'}^{q,z,a})$ in the problem \eqref{objectiveFunction}-\eqref{last_cons} to constants such that they satisfy the constraint in \eqref{integer_cons}. Note that this operation is equivalent to restricting the search space in \eqref{objective2}-\eqref{constraint2} to $k$-FSCs with a specific deterministic transition function, where the transition function satisfies $\delta({q'|q,z,a})$$=$$u(\delta_{q'}^{q,z,a})$. The resulting optimization problem has decision variables $\nu({\bf{s}})$, $\eta({\bf{s}})$, and $u(\gamma_a^{q,z})$, i.e., $u(\delta_{q'}^{q,z,a})$ is not a variable anymore. The resulting problem has bilinear constraints in \eqref{first_trilinear} and \eqref{reachabilityProb}, and hence, it is not a convex optimization problem. However, we can obtain a locally optimal solution to the resulting problem using a variant of the convex-concave-procedure (CCP) \cite{Yuille}. In particular, we employ \textit{penalty CCP} which is introduced in \cite{Lipp} and used in the context of pMCs in \cite{Cubuktepe:10.1007/978-3-030-01090-4_10}. The algorithm described below follows closely from the one proposed in \cite{Cubuktepe:10.1007/978-3-030-01090-4_10}.

Penalty CCP algorithm takes five inputs: a threshold constant $\epsilon$$>$$0$, an initial penalty constant $\tau_0$$>$$0$, a multiplication factor $\mu$$>$$1$, a maximum penalty constant $\tau_{\max}$, and  initial estimates $\hat{\nu}_0({\bf{s}})$, $\hat{\eta}_0({\bf{s}})$, and $\hat{u}_0(\gamma_a^{q,z})$ for the variables $\nu({\bf{s}})$, $\eta({\bf{s}})$, and $u(\gamma_a^{q,z})$, respectively. Moreover, for each iteration $k$$\in$$\mathbb{Z}_+$ of the algorithm, we recursively define $\tau_{k+1}$$:=$$\min\{\mu\tau_k,\tau_{\max}\}$. 

Let ${\bf{v}}$ denote an arbitrary tuple $({\bf{s}'},q,z,a)$$\in$$\mathcal{S}_{\mathbf{M},k}\times \mathcal{Q}\times \mathcal{Z}\times \mathcal{A}$. For each ${\bf{v}}$, we introduce two new variables $\Phi_{\nu,{\bf{v}}}$$\geq$$0$ and $\Phi_{\eta, {\bf{v}}}$$\geq$$0$. The introduced variables are typically referred to as \textit{slack variables} and quantify the infeasibility of the constraints in \eqref{first_trilinear}-\eqref{reachabilityProb} \cite{Lipp}. In particular, when $\sum_{\bf{v}}(\Phi_{\eta, {\bf{v}}}+\Phi_{\nu,{\bf{v}}})$$=$$0$, the output of the penalty CCP algorithm becomes feasible for the problem in \eqref{objectiveFunction}-\eqref{last_cons}.

At each iteration $k$$\in$$\mathbb{Z}_+$ of the algorithm, we first convexify the constraints in \eqref{first_trilinear}-\eqref{reachabilityProb} (explained below). We then solve the resulting convex optimization problem by replacing the objective function \eqref{objectiveFunction} with
\begin{align*}
     \underset{\nu,u, \eta}{\text{maximize}} \qquad  \nu({\bf{s}}_{I})-\tau_k\sum_{\bf{v}}(\Phi_{\eta, {\bf{v}}}+\Phi_{\nu,{\bf{v}}}).
\end{align*}
Intuitively, the second term in the above objective function is a penalty term which encourages the algorithm to output feasible solutions for the original problem in \eqref{objectiveFunction}-\eqref{last_cons}.

Let ${Val}_k$ be the optimal value of the problem described above. We terminate the algorithm if $\lvert {Val}_k-{Val}_{k-1}\rvert$$<$$\epsilon$ and the optimal solution satisfies $\sum_{\bf{v}}(\Phi_{\eta, {\bf{v}}}+\Phi_{\nu,{\bf{v}}})$$=$$0$. Otherwise, we set the optimal decision variables $\nu^{\star}({\bf{s}})$, $\eta^{\star}({\bf{s}})$, and $u^{\star}(\gamma_a^{q,z})$ for the current iteration as the estimates $\hat{\nu}_{k+1}({\bf{s}})$, $\hat{\eta}_{k+1}({\bf{s}})$, and $\hat{u}_{k+1}(\gamma_a^{q,z})$ for the successive iteration, and solve the resulting optimization problem. The procedure explained above has no theoretical convergence guarantees to a feasible solution \cite{Lipp}, i.e., a solution that satisfies $\sum_{\bf{v}}(\Phi_{\eta, {\bf{v}}}+\Phi_{\nu,{\bf{v}}})$$=$$0$. However, any feasible solution that is obtained through the above procedure is guaranteed to be locally optimal for the problem in \eqref{objectiveFunction}-\eqref{last_cons}. In practice, we observe that the penalty CCP usually converges to a feasible solution.

In what follows, we explain the convexification procedure for the constraint in \eqref{first_trilinear}; the convexification of \eqref{reachabilityProb} is performed by following the same procedure. Note that the last term on the right hand side of \eqref{first_trilinear} is the summation of bilinear terms $c(s,s',a,z,u) \nu({\bf{s'}}) u(\gamma_a^{q,z})$ where $c(s,s',a,z,u)$ is a constant such that
\begin{align*}
    c(s,s',a,z,u):={O}_{s,z} \,{P}_{s,a,s'} \, u(\delta_{q'}^{q,z,a}).
\end{align*}
 With an abuse of notation, we denote $c(s,s',a,z,u)$ by $c$. As explained in \cite{Cubuktepe:10.1007/978-3-030-01090-4_10}, a bilinear function $f(x,y)$$=$$2Cxy$, where $C$ is a constant, can be written as a difference of convex functions $f(x,y)$$=$$f_1(x,y)$$-$$f_2(x,y)$ where $f_1(x,y)$$=$$C(x$$+$$y)^2$ and $f_2(x,y)$$=$$C(x^2+y^2)$. Since we have a constraint of the form $0$$\leq$$L^u({\bf{s}})$$+$$f(x,y)$ in \eqref{first_trilinear}, we linearize the function $f_1(x,y)$ around the point $\hat{\nu}_k({\bf{s}})$ and $\hat{u}_k(\gamma_a^{q,z})$. Specifically, at iteration $k$$\in$$\mathbb{Z}_+$, we replace each bilinear term $c\nu({\bf{s'}}) u(\gamma_a^{q,z})$ in \eqref{first_trilinear} with 
\begin{align*}
&{\frac{c}{2}}\Big(\hat{\nu}_k({\bf{s}'})+\hat{u}_k(\gamma_a^{q,z})\Big)^2-\frac{c}{2}\Big((\nu({\bf{s}'}))^2+(u(\gamma_a^{q,z}))^2\Big)
\\
& +c\Big(\hat{\nu}_k({\bf{s}'})+\hat{u}_k(\gamma_a^{q,z})\Big)\Big(\nu({\bf{s}'})-\hat{\nu}_k({\bf{s}'})\Big) \\
&+ c\Big(\hat{\nu}_k({\bf{s}'})+\hat{u}_k(\gamma_a^{q,z})\Big)\Big(u(\gamma_a^{q,z})-\hat{u}_k(\gamma_a^{q,z}) \Big)+\Phi_{\nu,{\bf{v}}}.
\end{align*}
Note that the above expression is concave in the variables $\nu({\bf{s'}})$ and $u(\gamma_a^{q,z})$. Therefore, the problem resulting from the replacement of the bilinear terms with the above expression is a convex optimization problem.

\subsection{Finite-State Controller Synthesis: A Monotonocity Result}
In the previous section, we presented an algorithm to solve the problem in \eqref{objectiveFunction}-\eqref{last_cons}, which requires one to set the variables $u(\delta_{q'}^{q,z,a})$ to constants that satisfy the constraint in \eqref{integer_cons}. As discussed earlier, the choice of these constants establishes the memory transition structure of the $k$-FSCs over which the entropy maximization is performed. In this section, we present a particular memory transition function which has a monotonocity property. That is, when this memory transition function is used, by increasing the number of memory states, one can only increase the optimal value of the optimization problem in \eqref{objective2}-\eqref{constraint2}. 

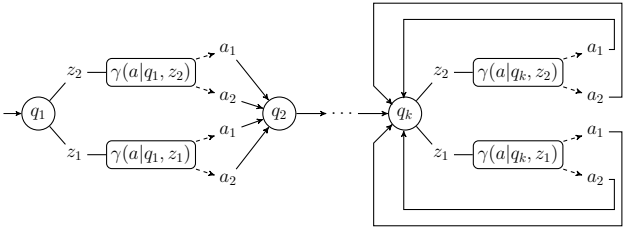
\begin{figure}[t!]\centering
\scalebox{0.5}{
\begin{tikzpicture}[->, >=stealth', auto, semithick, node distance=2cm]

     \tikzstyle{every state}=[fill=white,draw=black,thick,text=black,scale=1]

     \node[state,initial,initial text=] (q_1) {\Large{$q_1$}};
     \node        (z_1) [below right=5mm and 3mm of q_1] {\Large{$z_1$}};
     \node[draw=black, thick, rounded corners]   (gamma_1) [right=5mm of z_1] {\Large{$\gamma(a \lvert q_1,z_1)$}};
     \node   (a_11) [above right=-0.5 mm and 5mm of gamma_1] {\Large{$a_1$}};
     \node   (a_12) [below right=0 mm and 5mm of gamma_1] {\Large{$a_2$}};
     \node        (z_2) [above right=5mm and 3mm of q_1] {\Large{$z_2$}};
     \node[draw=black, thick, rounded corners]   (gamma_2) [right=5mm of z_2] {\Large{$\gamma(a \lvert q_1,z_2)$}};
     \node   (a_21) [above right=-0.5 mm and 5mm of gamma_2] {\Large{$a_1$}};
     \node   (a_22) [below right=0 mm and 5mm of gamma_2] {\Large{$a_2$}};
     \node[state] (q_2) [right=55mm of q_1]  {\Large{$q_2$}};
     \node[] (dots) [right=8mm of q_2]  {\Large{$\ldots$}};
      \node[state] (q_k) [right=8mm of dots]  {\Large{$q_k$}};
      
   \node        (z_k1) [below right=5mm and 3mm of q_k] {\Large{$z_1$}};
     \node[draw=black, thick, rounded corners]   (gamma_k1) [right=5mm of z_k1] {\Large{$\gamma(a \lvert q_k,z_1)$}};
     \node   (a_k11) [above right=-0.5 mm and 5mm of gamma_k1] {\Large{$a_1$}};
     \node   (a_k12) [below right=0 mm and 5mm of gamma_k1] {\Large{$a_2$}};
     \node   (z_k2) [above right=5mm and 3mm of q_k] {\Large{$z_2$}};
     \node[draw=black, thick, rounded corners]   (gamma_k2) [right=5mm of z_k2] {\Large{$\gamma(a \lvert q_k,z_2)$}};
     \node   (a_k21) [above right=-0.5 mm and 5mm of gamma_k2] {\Large{$a_1$}};
     \node   (a_k22) [below right=0 mm and 5mm of gamma_k2] {\Large{$a_2$}};

    \draw[thick,-] (q_1) -- (z_1);
    \draw[thick,-] (z_1) -- (gamma_1);
    \draw[thick,dashed, ->] (gamma_1) -- (a_11);
    \draw[thick,dashed, ->] (gamma_1) -- (a_12);
    \draw[thick,->] (a_11) -- (q_2);
    \draw[thick,->] (a_12) -- (q_2);
    \draw[thick,-] (q_1) -- (z_2);
    \draw[thick,-] (z_2) -- (gamma_2);
    \draw[thick,dashed, ->] (gamma_2) -- (a_21);
    \draw[thick,dashed, ->] (gamma_2) -- (a_22);
    \draw[thick,->] (a_21) -- (q_2);
    \draw[thick,->] (a_22) -- (q_2);
    \draw[thick,->] (q_2) -- (dots);
    \draw[thick,->] (dots) -- (q_k);
    
    \draw[thick,-] (q_k) -- (z_k1);
    \draw[thick,-] (z_k1) -- (gamma_k1);
    \draw[thick,dashed, ->] (gamma_k1) -- (a_k11);
    \draw[thick,dashed, ->] (gamma_k1) -- (a_k12);
    \draw[thick,-] (q_k) -- (z_k2);
    \draw[thick,-] (z_k2) -- (gamma_k2);
    \draw[thick,dashed, ->] (gamma_k2) -- (a_k21);
    \draw[thick,dashed, ->] (gamma_k2) -- (a_k22);
    \draw[thick,->] (a_k11) -| ++(0.7,-2.5) -| ++(-6.6,2.2) -- ++(0.5,0.5);
    \draw[thick,->] (a_k22) -| ++(0.7,2.5) -| ++(-6.6,-2.2) -- ++(0.5,-0.5);
     \draw[thick,->] (a_k12) -| ++(0.5,-0.8) -| ++(-5.6,2.1);
     \draw[thick,->] (a_k21) -| ++(0.5,0.8) -| ++(-5.6,-2.1);
\end{tikzpicture}}
\caption{ An illustration of the proposed deterministic $k$-FSC structure. Regardless of the received observations and taken actions, the controller transitions from the memory state $q_i$ to $q_{i+1}$ for all $i$$<$$k$, and finally, stays in the state $q_k$ indefinitely.  }\label{deterministic_controller_example}
\end{figure}

For a POMDP $\mathbf{M}$, consider a $k$-FSC $\mathbf{C}$$=$$(\mathcal{Q},q_1,\gamma,\overline{\delta})$ with the memory transition function $\overline{\delta}$$:$$\mathcal{Q}$$\times$$\mathcal{Z}$$\times$$\mathcal{A}$$\rightarrow$$\Delta(\mathcal{Q})$ such that
\begin{align}\label{LastLoopDetFSC}
\begin{cases}
        \overline{\delta}(q_{i+1}|q_{i},z,a)=1 & \forall z \in \mathcal{Z}, a \in \mathcal{A}, 1 \leq i < k \\
        \overline{\delta}(q_{k}|q_{k},z,a)=1  &\forall z \in \mathcal{Z}, a \in \mathcal{A}\\
        \overline{\delta}(q_i|q_j,z,a)=0 & \text{otherwise}.
        \end{cases}
\end{align}
We present an illustration of the $k$-FSC described above in Fig. \ref{deterministic_controller_example}. Intuitively, the transition function $\overline{\delta}$ represents a finite horizon memory. In the first $k$$-$$1$ steps, the agent summarizes the set $\mathcal{H}^i$ of system histories using the memory state $q_i$ and makes a decision based on the decision function $\gamma(a | q_i, z)$. For the rest of the process, the agent stays in the memory state $q_k$ and follows a memoryless strategy by making stationary decisions based on the decision function $\gamma(a | q_k, z)$.

Let $\overline{\mathcal{F}}_{k}(\mathbf{M})$$\subset$$\mathcal{F}_{k}^{det}(\mathbf{M})$ be the set of $k$-FSCs whose memory transition function is given in \eqref{LastLoopDetFSC}. Additionally, let 
 \begin{align*}
  E_{k,\max}:=\max_{\mathbf{C}\in\overline{\mathcal{F}}_{k}(\mathbf{M})}\ \ \sum_{t=2}^{\infty}\beta^{t-2}H^{\mathbf{C}}(S_{t} | S^{t-1}).
 \end{align*}

{\setlength{\parindent}{0cm}
{\begin{lemma}\label{lemma:pMCmonotonic}
For all $j$$\leq$$k$, we have $ E_{j,\max}$$\leq$$E_{k,\max}$.
\end{lemma}}}

The above monotonocity result establishes that, by \textit{fixing} the memory transition function to $\overline{\delta}$ given in \eqref{LastLoopDetFSC}, one can obtain nondecreasing maximum entropy values by increasing the number of memory states in the $k$-FSC. Such a monotonocity result hold due to the specific structure of the transition function $\overline{\delta}$ and may not hold if one considers transition functions other than $\overline{\delta}$. For example, consider a \textit{periodic} transition function $\widetilde{\delta}$ such that $\widetilde{\delta}(q_{j}|q_{i},z,a)$$=$$\overline{\delta}(q_{j}|q_{i},z,a)$ for all $1$$\leq$$i$$<$$k$, and $\widetilde{\delta}(q_{1}|q_{k},z,a)$$=$$1$. Due to the periodic structure of $\widetilde{\delta}$, the agent necessarily makes the same decisions in every $k$ steps. One can construct POMDP instances in which higher entropy values can be obtained by making the same decisions in every $k$ steps than in every $k$$+$$1$ steps. Therefore, the monotonocity result does not hold for the memory transition function $\widetilde{\delta}$.

Using the result of Lemma \ref{lemma:pMCmonotonic}, we can obtain a practical algorithm to synthesize an entropy-maximizing controller as follows. First, fix the number of memory states to an initial value $k$. Then, by setting $u(\delta_{q'}^{q,z,a})$$=$$\overline{\delta}(q'| q,z,a)$, find a local optimal solution to the problem in \eqref{objectiveFunction}-\eqref{last_cons}. Next, increase the number of memory states to $k$$+$$1$, solve the resulting problem, and compare the optimal value of the problem with the previous result. Repeat this procedure until the percent increase in the optimal value is below a predetermined threshold.

We note that since the algorithm described in the previous section computes a locally optimal solution for the problem in \eqref{objectiveFunction}-\eqref{last_cons}, the procedure described above has no theoretical guarantees to yield improving solutions for the problem in \eqref{objective2}-\eqref{constraint2}. However, we observe that, in practice, the described procedure works considerably well.

\section{Numerical Examples}\label{section:Numerical Examples}
We now provide several numerical examples to demonstrate the relation between the maximum entropy, the time horizon, the threshold on the total reward, and the number of memory states in the FSC. 
We use the MOSEK \cite{mosek} solver with the CVX \cite{cvx} interface to solve the convex optimization problems. To improve the approximation of exponential cone constraints, we use the CVXQUAD \cite{cvxquad} package.

\subsection{The Relation Between the Maximum Entropy and the Expected Total Reward Threshold}
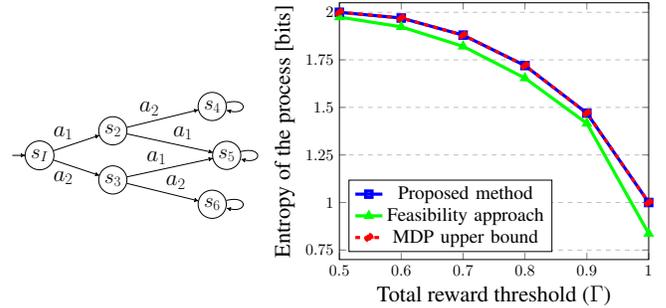
\begin{figure}[t!]
\begin{subfigure}{.2\textwidth}
    \centering
    \scalebox{0.35}{
        \begin{tikzpicture}[->, >=stealth', auto, semithick, node distance=2cm]
            \tikzstyle{every state}=[fill=white,draw=black,thick,text=black,scale=0.9]

            \node[state,initial,initial text=] (s_0) {\Huge{$s_I$}};
            \node[state] (s_1) [above right = 1.5mm and 20mm of s_0] {\Huge{$s_2$}};
            \node[state] (s_2) [below right = 1.5mm and 20mm of s_0] {\Huge{$s_3$}};
            \node[state] (s_3) [above right = 1.5mm and 30mm of s_1] {\Huge{$s_4$}};
            \node[state] (s_4) [right = 60mm of s_0] {\Huge{$s_5$}};
            \node[state] (s_5) [below right = 1.5mm and 30mm of s_2] {\Huge{$s_6$}};
        
            \path (s_0) edge node[above left]{\Huge{$a_1$}} (s_1)
                  (s_0) edge node[below left]{\Huge{$a_2$}} (s_2)
                  (s_1) edge node[above left]{\Huge{$a_2$}} (s_3)
                  (s_1) edge node[above right]{\Huge{$a_1$}} (s_4)
                  (s_2) edge node[above left]{\Huge{$a_1$}} (s_4)
                  (s_2) edge node[above right]{\Huge{$a_2$}} (s_5)
                  (s_3) edge[loop right] node[right]{} (s_3)
                  (s_4) edge[loop right] node[right]{} (s_4)
                  (s_5) edge[loop right] node[right]{} (s_5);
                  
        \end{tikzpicture}}
\end{subfigure}\hspace{-2mm}
\begin{subfigure}{.15\textwidth}
    \centering
    \scalebox{0.6}{
    \begin{tikzpicture}
        \begin{axis}[
            xlabel={\Large{Total reward threshold ($\Gamma$)}},
            ylabel={\Large{Entropy of the process [bits]}},
            xmin=0.5,xmax=1,
            ymin=0.70,ymax=2.05,
            xtick={0.5,0.6,0.7,0.8,0.9,1},
            ytick={0.75,1.00,1.25,1.50,1.75,2.00},
            legend pos=south west,
            ymajorgrids=true,
            grid style=dashed]
            \addplot[color=blue,mark=square,line width=2pt]
                coordinates{(0.5,2.00)(0.6,1.97)(0.7,1.88)(0.8,1.72)(0.9,1.47)(1.0,1.00)};
            \addplot[color=green,mark=triangle,line width=2pt]
                coordinates{(0.5,1.9750)(0.6,1.9235)(0.7,1.8209)(0.8,1.6540)(0.9,1.416)(1.0,0.8367)};
            \addplot[color=red,mark=*,dashed,line width=2pt]
                coordinates{(0.5,2.00)(0.6,1.97)(0.7,1.88)(0.8,1.72)(0.9,1.47)(1.0,1.00)};
            \legend{\large{Proposed method}, \large{Feasibility approach}, \large{MDP upper bound}}
        \end{axis}
    \end{tikzpicture}}
\end{subfigure}
    \caption{Example illustrating the relation between the maximum entropy and the threshold $\Gamma$. (Left) A POMDP instance with 6 states and a single observation. (Right) The entropy of the processes induced by the synthesized controllers as a function of the threshold $\Gamma$.}
    \label{fig:simpleExample}
\end{figure}

We first consider a POMDP instance with 6 states shown in Fig. \ref{fig:simpleExample} (Left). There is a single observation $\mathcal{Z}$$=$$\{z_1\}$, yielding $\mathcal{O}_{s,z_1}$$=$$1$ for all states \textit{s}. We use a deterministic 2-FSC whose memory transition function $\delta$ is given in \eqref{LastLoopDetFSC}. Since there is only one observation, the synthesized controller is an open-loop controller. We suppose that the agent aims to reach state $s_{4}$ and encode this objective by defining a reward function ${R}$ such that ${R}(s_2,a_1)$$=$${R}(s_3,a_1)$$=$$1$ and ${R}(s,a)$$=$$0$ otherwise.

We investigate the effect of the threshold $\Gamma$ in \eqref{constraint2} on the maximum entropy by synthesizing controllers for values between $\Gamma$$=$$0.5$ and $\Gamma$$=$$1$. For each value of $\Gamma$, we use the memory transition function in \eqref{LastLoopDetFSC}, run the optimization problem given in Section \ref{FSCsynthesis} 10 times by randomly initializing the CCP, and report the best result of  these 10 trials. For each $\Gamma$, we plot the maximum entropy of the stochastic process induced by the synthesized controllers in Fig. \ref{fig:simpleExample} (Right). For comparison, we synthesize controllers by solving a feasibility problem given in \cite{Cubuktepe:10.1007/978-3-030-01090-4_10}. We obtain the feasibility problem from \eqref{objectiveFunction}-\eqref{last_cons} by removing the entropy constraint in (\ref{first_trilinear}) and replacing the objective function in (\ref{objectiveFunction}) with a constant value.

The proposed approach yields the globally optimal controller by attaining a tight bound on $\Gamma$. The global optimality of the controller is evident in Fig. \ref{fig:simpleExample} (Right), as the entropy of the proposed approach exactly matches that of the underlying MDP for each value of $\Gamma$. Because the feasibility problem only seeks to find a feasible instantiation of the parameters that satisfy the expected total reward constraints in (\ref{reachabilityProb})-(\ref{initialStateProb}), the entropy of the induced stochastic processes is less than the maximum attainable entropy.

\subsection{The Relation Between the Maximum Entropy and the Number of Memory States}

\begin{figure}[t!]
\begin{subfigure}{\linewidth}
    \centering\scalebox{0.5}{
    \begin{tikzpicture}[->, >=stealth', auto, semithick, node distance=2cm]
        \tikzstyle{every state}=[fill=white,draw=black,thick,text=black,scale=0.8]

        \node[state,initial,initial text=] (s_2) {\Huge{$s_I$}};
        \node[state] (s_1) [above = 5mm of s_2] {\Huge{$s_2$}};
        \node[state] (s_3) [below = 5mm of s_2] {\Huge{$s_3$}};
        
        \node[state] (s_4) [right = 8mm of s_1] {\Huge{$s_4$}};
        \node[state] (s_7) [right = 8mm of s_4] {\Huge{$s_7$}};
        \node[state] (s_10) [right = 8mm of s_7] {\Huge{$s_{10}$}};
        \node[state] (s_13) [right = 8mm of s_10] {\Huge{$s_{13}$}};
        
        \node[state] (s_5) [right = 8mm of s_2] {\Huge{$s_5$}};
        \node[state] (s_8) [right = 8mm of s_5] {\Huge{$s_8$}};
        \node[state] (s_11) [right = 8mm of s_8] {\Huge{$s_{11}$}};
        \node[state] (s_14) [right = 8mm of s_11] {\Huge{$s_{14}$}};
        
        \node[state] (s_6) [right = 8mm of s_3] {\Huge{$s_6$}};
        \node[state] (s_9) [right = 8mm of s_6] {\Huge{$s_9$}};
        \node[state] (s_12) [right = 8mm of s_9] {\Huge{$s_{12}$}};
        \node[state] (s_15) [right = 8mm of s_12] {\Huge{$s_{15}$}};
        
        \node (dots_1) [right = 1mm of s_2] {$\cdots$};
        \node (dots_2) [right = 28mm of dots_1] {$\cdots$};
        \node (dots_3) [right = 13mm of dots_2] {$\cdots$};
        
        \node (dots_4) [right = 1mm of s_4] {$\cdots$};
        \node (dots_5) [right = 11mm of dots_4] {$\cdots$};
        \node (dots_6) [right = 12mm of dots_5] {$\cdots$};
        
        \node (dots_7) [right = 1mm of s_6] {$\cdots$};
        \node (dots_8) [right = 11mm of dots_7] {$\cdots$};
        \node (dots_9) [right = 12mm of dots_8] {$\cdots$};
        
        \path (s_13) edge[loop right] node[right]{} (s_13)
            (s_14) edge[loop right] node[right]{} (s_14)
            (s_15) edge[loop right] node[right]{} (s_15);
        
        \path (s_1) edge node[above]{\large{$a_1$}} (s_4)
            (s_1) edge node[above]{\large{$a_2$}} (s_5)
            (s_1) edge node[left]{\large{$a_3$}} (s_2);
            
        \path (s_5) edge node[above]{\large{$a_1$}} (s_7)
            (s_5) edge node[above]{\large{$a_2$}} (s_8)
            (s_5) edge node[above]{\large{$a_3$}} (s_9);
            
        \path (s_3) edge node[above]{\large{$a_1$}} (s_6)
            (s_3) edge node[above]{\large{$a_2$}} (s_5)
            (s_3) edge node[left]{\large{$a_3$}} (s_2);
             
        \end{tikzpicture}}
        \caption{A POMDP instance with 15 states and a single observation.}
\end{subfigure}
\vspace{0.5cm}

\begin{subfigure}{\linewidth}
    \centering\scalebox{0.5}{
        \begin{tikzpicture}[->, >=stealth', auto, semithick, node distance=2cm]
        \tikzstyle{every state}=[fill=white,draw=black,thick,text=black,scale=0.8]

        \node[state,initial,initial text=] (s_2) {\Huge{$s_I$}};
        \node[state] (s_1) [above = 5mm of s_2] {\Huge{$s_2$}};
        \node[state] (s_3) [below = 5mm of s_2] {\Huge{$s_3$}};
        
        \node[state] (s_4) [right = 8mm of s_1] {\Huge{$s_4$}};
        \node[state] (s_7) [right = 8mm of s_4] {\Huge{$s_7$}};
        \node[state] (s_10) [right = 8mm of s_7] {\Huge{$s_{10}$}};
        \node[state] (s_13) [right = 8mm of s_10] {\Huge{$s_{13}$}};
        
        \node[state] (s_5) [right = 8mm of s_2] {\Huge{$s_5$}};
        \node[state] (s_8) [right = 8mm of s_5] {\Huge{$s_8$}};
        \node[state] (s_11) [right = 8mm of s_8] {\Huge{$s_{11}$}};
        \node[state] (s_14) [right = 8mm of s_11] {\Huge{$s_{14}$}};
        
        \node[state] (s_6) [right = 8mm of s_3] {\Huge{$s_6$}};
        \node[state] (s_9) [right = 8mm of s_6] {\Huge{$s_9$}};
        \node[state] (s_12) [right = 8mm of s_9] {\Huge{$s_{12}$}};
        \node[state] (s_15) [right = 8mm of s_12] {\Huge{$s_{15}$}};
        
        \draw[line width=0.3mm, red] (s_2) -- (s_4);
        \draw[line width=0.3mm, red] (s_2) -- (s_5);
        \draw[line width=0.3mm, red] (s_2) -- (s_6);
        
        \draw[line width=0.3mm, red] (s_4) -- (s_5);
        \draw[line width=0.3mm, red] (s_4) -- (s_7);
        \draw[line width=0.3mm, red] (s_4) -- (s_8);
        
        \draw[line width=0.3mm, red] (s_5) -- (s_7);
        \draw[line width=0.3mm, red] (s_5) -- (s_8);
        \draw[line width=0.3mm, red] (s_5) -- (s_9);
        
        \draw[line width=0.3mm, red] (s_6) -- (s_5);
        \draw[line width=0.3mm, red] (s_6) -- (s_8);
        \draw[line width=0.3mm, red] (s_6) -- (s_9);
        
        \draw[line width=0.3mm, red] (s_7) -- (s_8);
        \draw[line width=0.3mm, red] (s_7) -- (s_10);
        \draw[line width=0.3mm, red] (s_7) -- (s_11);
        
        \draw[line width=0.3mm, red] (s_8) -- (s_10);
        \draw[line width=0.3mm, red] (s_8) -- (s_11);
        \draw[line width=0.3mm, red] (s_8) -- (s_12);
        
        \draw[line width=0.3mm, red] (s_9) -- (s_8);
        \draw[line width=0.3mm, red] (s_9) -- (s_11);
        \draw[line width=0.3mm, red] (s_9) -- (s_12);
        
        \draw[line width=0.60mm, red] (s_10) -- (s_14);
        \draw[line width=0.60mm, red] (s_10) -- (s_11);
        
        \draw[line width=0.9mm, red] (s_11) -- (s_14);
        
        \draw[line width=0.9mm, red] (s_12) -- (s_14);
        \draw[line width=0.60mm, red] (s_12) -- (s_11);
        
        \path[line width=0.9mm, red] (s_14) edge [out=45,in=315,looseness=6] node[above] {} (s_14);
        \end{tikzpicture}}
        \caption{Trajectories under the synthesized entropy maximizing 6-FSC. Edge thicknesses indicate the transition probabilities.}
\end{subfigure}
    \caption{Example illustrating the relation between the maximum entropy and the number of memory states. }
    \label{fig:POMDPoptController}
\end{figure}
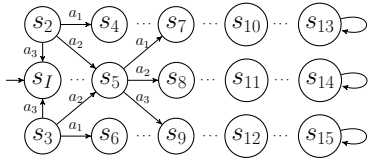
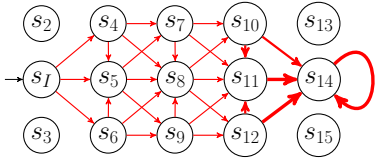

We now consider a POMDP instance with 15 states shown in Fig. \ref{fig:POMDPoptController} (a). As in the previous example, there is only a single observation $\mathcal{Z}$$=$$\{z_{1}\}$ yielding $\mathcal{O}_{s,z_{1}}$$=$$1$ for all states $s$. We suppose that the agent aims to reach $s_{14}$ with probability 1. To encode this objective, we set $\Gamma$$=$$1$ with ${R}(s_{10},a_{2})$$=$$1$, ${R}(s_{11},a_{2})$$=$$1$, ${R}(s_{12},a_{2})$$=$$1$, and ${R}(s,a)$$=$$0$ otherwise.

\begin{figure}[t]
    \centering
    \scalebox{0.6}{
    \begin{tikzpicture}
        \begin{axis}[
            xlabel={\Large{Number of memory states ($k$)}},
            ylabel={\Large{Entropy of the process [bits]}},
            xmin=0.95,xmax=6.05,
            ymin=-0.05,ymax=7.25,
            xtick={0,1,2,3,4,5,6},
            ytick={0,1.0,2.0,3.0,4.0,5.0,6.0,7.0},
            ymajorgrids=true,
            grid style=dashed,
            legend pos=south east]
            \addplot[color=blue,mark=square,line width=2pt]
                coordinates{(1,0)(2,1.59)(3,3.17)(4,4.76)(5,5.17)(6,5.17)};
            \addplot[color=red,mark=none,line width=2pt]
                coordinates{(0,6.89)(2,6.89)(3,6.89)(4,6.89)(5,6.89)(7,6.89)};
            \legend{\large{Proposed method}, \large{MDP upper bound}};
        \end{axis}
    \end{tikzpicture}}
    \caption{The maximum entropy of the stochastic processes induced by deterministic $k$-FSCs. The maximum entropy is monotonically nondecreasing in the number of memory states.}
    \label{fig:EntropyVsk}
\end{figure}
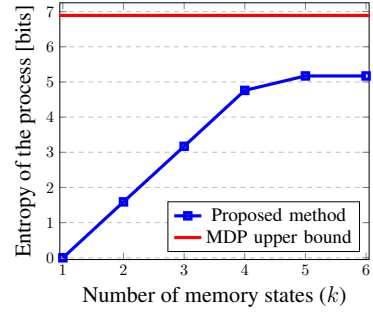

We study the relation between the number of memory states and the maximum entropy of the induced pMC by synthesizing controllers with $k$$=$$1,\ldots,6$ memory states. We run the optimization problem given in Section \ref{FSCsynthesis} 10 times by randomly initializing the CCP and report the best result of these 10 trials. Fig. \ref{fig:EntropyVsk} shows the maximum entropy of the induced stochastic process for each $k$. Moreover, in Fig. \ref{fig:POMDPoptController} (b), we illustrate the entropy-maximizing 6-FSC, where edge weights correspond to the probability of action selection. 

As can be seen in Fig. \ref{fig:EntropyVsk}, the maximum entropy of the POMDP under 1-FSCs is 0 since, under such controllers, the agent needs to follow the shortest path to the goal for collecting a total reward of 1. When we increase the number of memory states, the agent starts to randomize its action selection. However, after 5 memory states, memory does not affect the maximum entropy of the induced process. Note in Fig. \ref{fig:EntropyVsk} that, unlike the previous example, a gap between the maximum entropy of the MDP and that of the induced pMC remains. The maximum entropy of the POMDP must lie within this gap. Finally, this example demonstrates that the maximum entropy of the induced process is monotonically nondecreasing in the number of memory states when we consider controllers with memory transition function given in \eqref{LastLoopDetFSC}.

\subsection{The Relation Between the Maximum Entropy and the Time Horizon}\label{example_3}

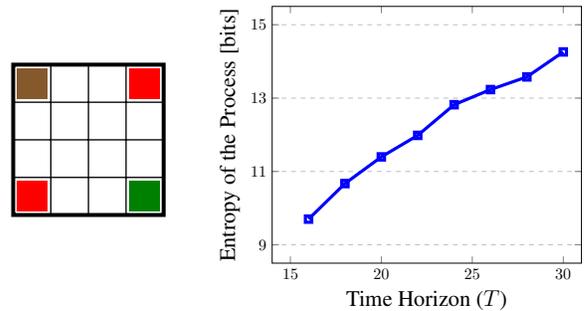
\begin{figure}[b!]
\begin{subfigure}[b]{0.17\textwidth}
\centering\raisebox{1.3cm}{
\scalebox{0.5}{
\begin{tikzpicture}
\draw[black,line width=1pt] (0,0) grid[step=1] (4,4);
\draw[black,line width=3pt] (0,0) rectangle (4,4);

\initialstate{0}{3}  \goalstate{3}{0}
\StaticObstacle{3}{3} \StaticObstacle{0}{0}
\end{tikzpicture}
}}
\end{subfigure}
\begin{subfigure}[b]{0.2\textwidth}
\scalebox{0.6}{
    \begin{tikzpicture}
        \begin{axis}[
            xlabel={\Large{Time Horizon ($T$)}},
            ylabel={\Large{Entropy of the Process [bits]}},
            xmin=14,xmax=31,
            ymin=8.5,ymax=15.5,
            xtick={15,20,25,30},
            ytick={9,11,13,15},
            ymajorgrids=true,
            grid style=dashed,
            legend pos=south east]
            \addplot[color=blue,mark=square,line width=2pt]
                coordinates{(16,9.7008)(18,10.6696)(20,11.3973)(22,11.9841)(24,12.8198)(26,13.2323)(28,13.5781)(30,14.2590)};
        \end{axis}
    \end{tikzpicture}}
\end{subfigure}
\caption{Example illustrating the relation between the maximum entropy and the time horizon. (Left) Grid world with 16 states and 9 observations. The agent starts from the brown state and aims to reach the green state while avoiding red states. (Right) The entropy of the processes induced by the synthesized $1$-FSCs as a function of the time horizon. }
\label{simulation_tikz_4_by_4}
\end{figure}

In this example, we examine the relation between the maximum entropy of a POMDP and the number of time steps in the finite-horizon problem. We consider the POMDP whose state-space is given by the $4$$\times$$4$ grid world shown in Fig. \ref{simulation_tikz_4_by_4} (Left). The brown state is the unique initial state of the agent, the red states are error states to be avoided, and the green state is the target state. In each state, the agent can select one of four possible actions: move left, move right, move up, or move down. Under the selected action, the agent transitions to its intended state with probability $0.95$$-$$\nicefrac{\epsilon}{3}$, slips to the left or to the right with probability $0.025$$-$$\nicefrac{\epsilon}{3}$, and slips backwards with probability $\epsilon$$=$$0.005$. If the agent were to transition off the grid world, it instead remains in its current state. 

The agent can make 9 possible observations: \emph{error (target) state to the left, error (target) state to the right, error (target) state above, error (target) state below}, and \emph{no observation}. The observations in each state are deterministic. For example, in the state to the left of the target state, the agent makes the observation \emph{target state to the right} with probability 1.

            

            

We investigate how the variation of the finite time horizon $T$ affects the maximum entropy of the POMDP. We use values of $T$ ranging from $T$$=$$16$ to $T$$=$$30$ time steps. The minimum expected total reward threshold is set to $\Gamma$$=$$0.9$, which we encode by setting the values of the reward function ${R}$ as ${R}((3,1),\text{right})$$=$${R}((4,2),\text{down})$$=$$0.95$$-$$\nicefrac{\epsilon}{3}$, ${R}((3,1),\text{up})$$=$${R}((4,2),\text{left})$$=$$0.025$$-$$\nicefrac{\epsilon}{3}$, and ${R}((3,1)$$,$$\text{left})$ =${R}((4,2),\text{up})$$=$$\epsilon$, and $0$ otherwise. We use a deterministic $1$-FSC and run the optimization problem given in Section \ref{FSCsynthesis} a total of 5 times for each value of $T$. In Fig. \ref{simulation_tikz_4_by_4} (Right), we plot the maximum entropy of the induced stochastic processes as a function of the finite time horizon $T$.

From Fig. \ref{simulation_tikz_4_by_4} (Right), we see that the maximum entropy of the POMDP increases with the size of the finite time horizon. Intuitively, a longer time horizon allows the agent to more uniformly distribute its actions. With a shorter time horizon, the agent must allocate more probability mass towards selecting actions that lead it down and to the right in order to reach the target state within the time horizon. For longer time horizons, the agent is able to more uniformly distribute its actions, yielding a higher maximum entropy.

\subsection{The Relation Between the Maximum Entropy and the Finite-Horizon Expected Total Reward}\label{example_4}
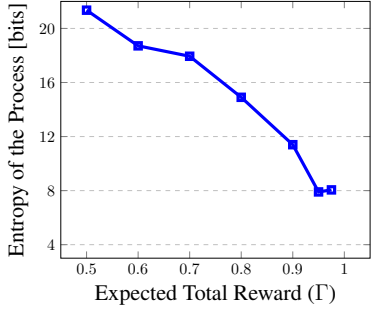
\begin{figure}[t!]
    \centering
    \scalebox{0.6}{
    \begin{tikzpicture}
        \begin{axis}[
            xlabel={\Large{Expected Total Reward ($\Gamma$)}},
            ylabel={\Large{Entropy of the Process [bits]}},
            xmin=0.45,xmax=1.05,
            ymin=3,ymax=22,
            xtick={0.4,0.5,0.6,0.7,0.8,0.9,1.0},
            ytick={4,8,12,16,20},
            ymajorgrids=true,
            grid style=dashed,
            legend pos=south east]
            \addplot[color=blue,mark=square,line width=2pt]
                coordinates{(0.5,21.3515)(0.6,18.7123)(0.7,17.9418)(0.8,14.9038)(0.9,11.3973)(0.95,7.9036)(0.975,8.0533)};
        \end{axis}
    \end{tikzpicture}}
    \caption{Relation between the maximum entropy of the induced stochastic process and the expected total reward $\Gamma$ for finite time horizon $T$$=$$20$.}
    \label{fig:EntropyVsGamma}
\end{figure}

In this example, we again consider the $4$$\times$$4$ grid world shown in Fig. \ref{simulation_tikz_4_by_4} (Left) with identical transition, observation, and reward functions as those used in Section \ref{example_3}. We now suppose that the agent must collect an expected total reward above some minimum threshold $\Gamma$ within a finite time horizon of $T$$=$$16$ time steps. We use values of $\Gamma$ varying from $\Gamma$$=$$0.5$ to $\Gamma$$=$$0.975$. Using a deterministic 1-FSC, we run the optimization problem given in Section \ref{FSCsynthesis} a total of 5 times for each value of $\Gamma$ and store the largest value of the maximum entropy observed. Fig. \ref{fig:EntropyVsGamma} plots the relation between the maximum entropy of the induced stochastic process as a function of the lower bound $\Gamma$ on the total expected reward.

Lower values of $\Gamma$ allow the agent to more uniformly distribute its probabilities for action selection, as the agent need not reach the blue state with high probability. For increasingly large values of $\Gamma$, the agent must select its actions such that it drives itself towards the blue state, reducing the maximum entropy of the induced stochastic process. As the value of $\Gamma$ approaches 1, the maximum entropy of the induced stochastic process begins to level off. For large values of $\Gamma$, and for $\Gamma$$=$$0.95$ and $\Gamma$$=$$0.975$ in particular, the synthesized controllers only slightly deviate from one another. Because the synthesized policies are nearly identical, the resulting maximum entropies of their respective induced stochastic processes are likewise nearly identical.

\subsection{Reaching a Target While Minimizing Predictability}\label{example_5}
We consider an agent that aims to reach a target in an adversarial environment. We model the environment as a $10$$\times$$10$ grid world, shown in Fig. 9 (Left), that consists of 4 \textit{rooms} and 4 \textit{doors} using which the agent can transition between the rooms. The rooms are numbered clockwise starting from the bottom left corner, and the doors are numbered clockwise starting from the door between room 1 and room 4. Finally, the thick black lines represent the walls in the environment. 

The agent observes its current room and its relative position with respect to the doors (36 total observations, 9 in each room). We illustrate the partition of a room with respect to the agent's observation function in Fig. 9 (Right). For example, if the agent is at the bottom left corner of the environment, its observation is \textit{room 1, below door 1, left of door 2}.

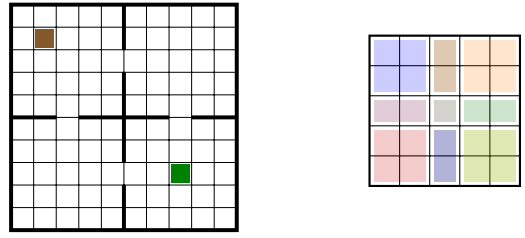
\begin{figure}[t!]
\centering
\begin{subfigure}[b]{0.25\textwidth}
\centering
\scalebox{0.3}{
\begin{tikzpicture}
\draw[black,line width=1pt] (0,0) grid[step=1] (10,10);
\draw[black,line width=5pt] (0,0) rectangle (10,10);
\draw[black,line width=5pt] (5,10) -- (5,8);
\draw[black,line width=5pt] (5,7) -- (5,3);
\draw[black,line width=5pt] (5,2) -- (5,0);
\draw[black,line width=5pt] (0,5) -- (2,5);
\draw[black,line width=5pt] (3,5) -- (7,5);
\draw[black,line width=5pt] (8,5) -- (10,5);
\initialstate{1}{8}  
\goalstate{7}{2}
\end{tikzpicture}
}
\end{subfigure}
\begin{subfigure}[b]{0.2\textwidth}
\centering\raisebox{0.6cm}{
\scalebox{0.4}{
\begin{tikzpicture}
\draw[black,line width=1pt] (0,0) grid[step=1] (5,5);
\draw[black,line width=2pt] (0,0) rectangle (5,5);
\obsa{0}{3}
\obsb{0}{0}
\obsc{3}{3}
\obsd{3}{0}
\obsf{0}{2}
\obsg{3}{2}
\obsq{2}{2}
\obsqq{2}{0}
\obsqqq{2}{3}
\end{tikzpicture}
}}
\end{subfigure}
\caption{ Motion planning example (Left) Grid world with 100 states and 36 observations. The agent starts from the brown state and aims to reach the green state. (Right) The partition of a room with respect to the agent's observation function. }

\label{simulation_tikz_10_by_10}
\end{figure}

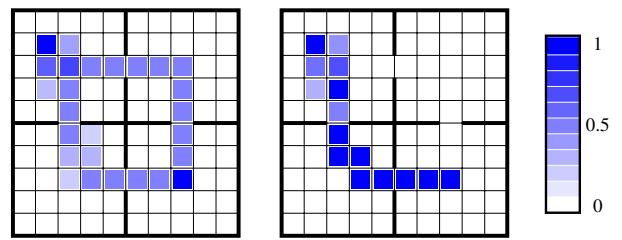
\begin{figure}[b!]
\centering
\begin{subfigure}[b]{0.21\textwidth}

\centering
\scalebox{0.3}{
\begin{tikzpicture}
\draw[black,line width=1pt] (0,0) grid[step=1] (10,10);
\draw[black,line width=5pt] (0,0) rectangle (10,10);
\draw[black,line width=5pt] (5,10) -- (5,8);
\draw[black,line width=5pt] (5,7) -- (5,3);
\draw[black,line width=5pt] (5,2) -- (5,0);
\draw[black,line width=5pt] (0,5) -- (2,5);
\draw[black,line width=5pt] (3,5) -- (7,5);
\draw[black,line width=5pt] (8,5) -- (10,5);

\fill[blue!100!white] (1+0.1,8+0.1) rectangle (1+0.9,8+0.9);
\fill[blue!37!white] (2+0.1,8+0.1) rectangle (2+0.9,8+0.9);

\fill[blue!63!white] (1+0.1,7+0.1) rectangle (1+0.9,7+0.9);
\fill[blue!73!white] (2+0.1,7+0.1) rectangle (2+0.9,7+0.9);
\fill[blue!50!white] (3+0.1,7+0.1) rectangle (3+0.9,7+0.9);
\fill[blue!50!white] (4+0.1,7+0.1) rectangle (4+0.9,7+0.9);
\fill[blue!50!white] (5+0.1,7+0.1) rectangle (5+0.9,7+0.9);
\fill[blue!50!white] (6+0.1,7+0.1) rectangle (6+0.9,7+0.9);
\fill[blue!50!white] (7+0.1,7+0.1) rectangle (7+0.9,7+0.9);

\fill[blue!27!white] (1+0.1,6+0.1) rectangle (1+0.9,6+0.9);
\fill[blue!49!white] (2+0.1,6+0.1) rectangle (2+0.9,6+0.9);
\fill[blue!50!white] (7+0.1,6+0.1) rectangle (7+0.9,6+0.9);

\fill[blue!50!white] (2+0.1,5+0.1) rectangle (2+0.9,5+0.9);
\fill[blue!50!white] (7+0.1,5+0.1) rectangle (7+0.9,5+0.9);

\fill[blue!50!white] (2+0.1,4+0.1) rectangle (2+0.9,4+0.9);
\fill[blue!18!white] (3+0.1,4+0.1) rectangle (3+0.9,4+0.9);
\fill[blue!50!white] (7+0.1,4+0.1) rectangle (7+0.9,4+0.9);

\fill[blue!31!white] (2+0.1,3+0.1) rectangle (2+0.9,3+0.9);
\fill[blue!30!white] (3+0.1,3+0.1) rectangle (3+0.9,3+0.9);
\fill[blue!50!white] (7+0.1,3+0.1) rectangle (7+0.9,3+0.9);

\fill[blue!20!white] (2+0.1,2+0.1) rectangle (2+0.9,2+0.9);
\fill[blue!50!white] (3+0.1,2+0.1) rectangle (3+0.9,2+0.9);
\fill[blue!50!white] (4+0.1,2+0.1) rectangle (4+0.9,2+0.9);
\fill[blue!50!white] (5+0.1,2+0.1) rectangle (5+0.9,2+0.9);
\fill[blue!50!white] (6+0.1,2+0.1) rectangle (6+0.9,2+0.9);
\fill[blue!100!white] (7+0.1,2+0.1) rectangle (7+0.9,2+0.9);

\end{tikzpicture}}
\end{subfigure}
\begin{subfigure}[b]{0.16\textwidth}
\centering{
\scalebox{0.3}{
\begin{tikzpicture}
\draw[black,line width=1pt] (0,0) grid[step=1] (10,10);
\draw[black,line width=5pt] (0,0) rectangle (10,10);
\draw[black,line width=5pt] (5,10) -- (5,8);
\draw[black,line width=5pt] (5,7) -- (5,3);
\draw[black,line width=5pt] (5,2) -- (5,0);
\draw[black,line width=5pt] (0,5) -- (2,5);
\draw[black,line width=5pt] (3,5) -- (7,5);
\draw[black,line width=5pt] (8,5) -- (10,5);

\fill[blue!100!white] (1+0.1,8+0.1) rectangle (1+0.9,8+0.9);
\fill[blue!43!white] (2+0.1,8+0.1) rectangle (2+0.9,8+0.9);

\fill[blue!56!white] (1+0.1,7+0.1) rectangle (1+0.9,7+0.9);
\fill[blue!69!white] (2+0.1,7+0.1) rectangle (2+0.9,7+0.9);

\fill[blue!31!white] (1+0.1,6+0.1) rectangle (1+0.9,6+0.9);
\fill[blue!100!white] (2+0.1,6+0.1) rectangle (2+0.9,6+0.9);

\fill[blue!50!white] (2+0.1,5+0.1) rectangle (2+0.9,5+0.9);

\fill[blue!100!white] (2+0.1,4+0.1) rectangle (2+0.9,4+0.9);

\fill[blue!100!white] (2+0.1,3+0.1) rectangle (2+0.9,3+0.9);
\fill[blue!100!white] (3+0.1,3+0.1) rectangle (3+0.9,3+0.9);

\fill[blue!100!white] (3+0.1,2+0.1) rectangle (3+0.9,2+0.9);

\fill[blue!100!white] (3+0.1,2+0.1) rectangle (3+0.9,2+0.9);
\fill[blue!100!white] (4+0.1,2+0.1) rectangle (4+0.9,2+0.9);
\fill[blue!100!white] (5+0.1,2+0.1) rectangle (5+0.9,2+0.9);
\fill[blue!100!white] (6+0.1,2+0.1) rectangle (6+0.9,2+0.9);
\fill[blue!100!white] (7+0.1,2+0.1) rectangle (7+0.9,2+0.9);

\end{tikzpicture}
}}
\end{subfigure}
\begin{subfigure}[b]{0.1\textwidth}
\centering\raisebox{0.3cm}{
\scalebox{0.3}{
\begin{tikzpicture}
  \colorlet{color min rgb}[rgb]{white}
  \colorlet{color max rgb}[rgb]{blue}
  \def\min{0}
  \def\max{100}

\fill[fill=blue!0!white,draw=white] (1cm,1 cm) rectangle +(1.5cm,0.8cm);
\fill[fill=blue!10!white,draw=white] (1cm,1.7 cm) rectangle +(1.5cm,0.8cm);
\fill[fill=blue!20!white,draw=white] (1cm,2.4 cm) rectangle +(1.5cm,0.8cm);
\fill[fill=blue!30!white,draw=white] (1cm,3.1 cm) rectangle +(1.5cm,0.8cm);
\fill[fill=blue!40!white,draw=white] (1cm,3.8 cm) rectangle +(1.5cm,0.8cm);
\fill[fill=blue!50!white,draw=white] (1cm,4.5 cm) rectangle +(1.5cm,0.8cm);
\fill[fill=blue!60!white,draw=white] (1cm,5.2 cm) rectangle +(1.5cm,0.8cm);
\fill[fill=blue!70!white,draw=white] (1cm,5.9 cm) rectangle +(1.5cm,0.8cm);
\fill[fill=blue!80!white,draw=white] (1cm,6.6 cm) rectangle +(1.5cm,0.8cm);
\fill[fill=blue!90!white,draw=white] (1cm,7.3 cm) rectangle +(1.5cm,0.8cm);
\fill[fill=blue!100!white,draw=white] (1cm,8 cm) rectangle +(1.5cm,0.8cm);

\draw[draw=black, line width = 3 pt] (1cm,1 cm) rectangle +(1.5cm,7.85cm);
\node at (3.3cm, 8 cm + 5mm) {\Huge 1};
\node at (3.3cm, 4.4 cm + 5mm) {\Huge 0.5};
\node at (3.3cm, 0.8 cm + 5mm) {\Huge 0};

\end{tikzpicture}}}
\end{subfigure}

\caption{ The expected number of times the agent visits each state under the synthesized controllers. (Left) Entropy-maximizing controller. (Right) Feasibility approach. }

\label{simulation_tikz_10_by_10_result}
\end{figure}

In Fig. 9 (Left), the brown state is the initial state of the agent, and the green state is the target state. We set the discount factor to $\beta$$=$$0.9$, and the expected total reward threshold to $\Gamma$$=$$\beta^{12}$, i.e., the agent needs to reach the target state in at most 12 steps. Note that 12 steps is the minimum number of steps to reach the target from the initial state. Hence, the agent can follow only the shortest trajectories to the target. 

We focus on $1$-FSCs and synthesize two controllers for the agent. We synthesize the first controller using the proposed approach based on the convex-concave procedure. For comparison, we also synthesize a controller by solving a feasibility problem given in \cite{Cubuktepe:10.1007/978-3-030-01090-4_10}. In Fig. 10, we demonstrate the expected number of times the agent visits each state under the synthesized controllers. 

As can be seen from Fig. 10, under the entropy-maximizing controller, the agent reaches the target state by passing through room 1 and room 3 with equal probability. Hence, it minimizes the predictability of the room it visits to an outside observer by maximizing the entropy of its trajectories. On the other hand, under the controller synthesized by the feasibility approach, the agent always reaches the target by passing through room 1. Consequently, it becomes trivial for an outside observer to predict the agent's trajectory.

\section{Conclusions}\label{section:Conclusions and Future Extensions}
We studied the synthesis of a controller which, from a given POMDP, induces a stochastic process with maximum entropy among the ones whose realizations accumulate a certain level of expected reward. Since the entropy maximization objective is considerably different than the traditionally used expected reward maximization objective, we first showed that the maximum entropy of a POMDP is upper bounded by the maximum entropy of its corresponding fully observable counterpart. Then, by restricting our attention to FSCs with deterministic memory transitions, we recast the entropy maximization problem as a so-called parameter synthesis problem for pMCs. We present a nonlinear optimization problem for the synthesis of an FSC that maximizes the entropy of a POMDP over all FSCs with the same number of memory states and deterministic memory transitions. Considering the intractability of finding a global optimal solution to the presented optimization problem, we proposed a convex-concave procedure approach to obtain a local optimal solution after setting the memory transition of FSCs to a fixed structure.

Even though finding a solution to the constrained entropy optimization problem is at least PSPACE-hard due to expected reward constraints, the computational complexity of the unconstrained entropy maximization problem is still an open problem. Additionally, developing the computational methods to synthesize a controller that maximizes the entropy of a POMDP over all FSCs with the same number of memory states may be a fruitful research direction. 

\bibliographystyle{IEEEtran}
\bibliography{ref}
\appendices
\section{} \label{section:Appendix}
In this appendix, we provide proofs for all theoretical results presented in the paper.

\textbf{Proof of Lemma \ref{writeValueFunction}.} For any $t$$<$$N$, we have
\begin{subequations}
\allowdisplaybreaks
    \begin{flalign*}
       {V}_{t,N}^{\pi}(s^{t}) = &H^{\pi}(S_{t+1}|S^{t}=s^{t}) + \sum_{k=t+1}^{N-1} H^{\pi}(S_{k+1}|S_{t}^{k},S^{t}=s^{t}) 
    \end{flalign*}
\end{subequations}
since the random variable $S_{t}$ is a part of the sequence $S^t$. The last term in the above equation satisfies
\begin{subequations}
\begin{flalign*}
    H^{\pi}(S_{k+1}|S_{t}^{k},S^{t}=s^{t}) &= H^{\pi}(S_{k+1}|S_{t+1}^{k},S^{t}=s^{t})\\ &=H^{\pi}(S_{k+1}|S_{t+2}^{k}, S_{t+1},S^{t}=s^{t})
\end{flalign*}
\end{subequations}
since $S_{t}^k$$=$$(S_t,S_{t+1},S^k_{t+2})$. Moreover, the term on the left hand side of the above equality satisfies
\begin{subequations}
\begin{flalign*}
 H^{\pi}(S_{k+1}|S_{t+2}^{k}, S_{t+1},S^{t}=s^{t}) = H^{\pi}(S_{k+1}|S_{t}^{k}, S_{t+1},S^{t}=s^{t})
\end{flalign*}
\end{subequations}
since the introduced conditioning on the sequence $(S_t,S_{t+1})$ does not change entropy as the random variable $S_{t+1}$ is already included in the conditioning and the value of the random variable $S_t$ is already fixed to $s_t$. Using the definition of conditional entropy \cite[Chapter 2]{Cover}, we obtain 
\begin{subequations}
\begin{flalign*}
 &H^{\pi}(S_{k+1}|S_{t}^{k}, S_{t+1},S^{t}=s^{t}) =\\ 
 &\sum_{s_{t+1} \in \mathcal{S}}  Pr(S_{t+1}=s_{t+1}|S^t=s^{t})H^{\pi}(S_{k+1}|S_{t}^{k},S^{t+1}=s^{t+1}).
\end{flalign*}
\end{subequations}
Note that, under the policy $\pi$, $Pr(S_{t+1}$$=$$s_{t+1}|S^t$$=$$s^{t})$ is equal to the realization probability $Pr^{\pi}(s^{t+1}|s^t)$. Furthermore, $Pr^{\pi}(s^{t+1}|s^t)$$>$$0$ for a given state history $s^{t+1}$$\in$$\mathcal{SH}^{t+1}$ if and only if $s^{t+1}$$=$$(s^t,s_{t+1})$ where $s_{t+1}$$\in$$\mathcal{S}$. As a result, we have 
\begin{subequations}
    \begin{flalign}
      & {V}_{t,N}^{\pi}(s^{t}) = H^{\pi}(S_{t+1}|S^{t}=s^{t})\nonumber &&\\
        & + \sum_{k=t+1}^{N-1}\sum_{s^{t+1} \in \mathcal{S}\mathcal{H}^{t+1}}  Pr^{\pi}(s^{t+1}|s^{t})H^{\pi}(S_{k+1}|S_{t}^{k},S^{t+1}=s^{t+1}) \label{SecondInd-e}&&\\
    & = H^{\pi}(S_{t+1}|S^{t}=s^{t})\nonumber &&\\
        & +\sum_{s^{t+1} \in \mathcal{S}\mathcal{H}^{t+1}}  Pr^{\pi}(s^{t+1}|s^{t})\sum_{k=t+1}^{N-1}H^{\pi}(S_{k+1}|S_{t}^{k},S^{t+1}=s^{t+1}) \label{SecondInd-f}&&\\
        & = H^{\pi}(S_{t+1}|S^{t}=s^{t}) \nonumber&& \\
        & + \sum_{s^{t+1} \in \mathcal{S} \mathcal{H}^{t+1}} Pr^{\pi}(s^{t+1}|s^{t}) {V}_{t+1,N}^{\pi}(s^{t+1})&& \label{SecondInd-g} \raisetag{20pt}
    \end{flalign}
\end{subequations}
where (\ref{SecondInd-f}) follows from the fact that the expression $Pr^{\pi}(s^{t+1}|s^{t})$ does not depend on $k$, and (\ref{SecondInd-g}) follows from the definition of the value function ${V}_{t,N}^{\pi}(s^{t})$: $\Box$

\textbf{Proof of Theorem \ref{thm:POMDPbounded}:} Recall that, for any $\pi$$\in$$\Pi(\mathbf{M})$, we have $V_{1,N}^{\pi}(s_I)$$=$$H^{\pi}(S_1,S_2,\ldots,S_N)$. We prove the claim by showing that 
\begin{align*}
    \sup_{\pi\in\Pi(\mathbf{M})} V_{1,N}^{\pi}(s_I)\leq  \sup_{\pi\in\Pi(\mathbf{M}_{fo})} V_{1,N}^{\pi}(s_I).
\end{align*}
We do so by considering the value function $ V_{t,N}^{\pi}(s^t)$ and performing induction on $t$. 

Denote the value function for an \textit{arbitrary} controller $\pi$$\in$$\Pi(\mathbf{M})$ as ${V}_{t,N}^{\pi}(s^{t})$ and the value function for the corresponding controller $\pi'$$=$$(\mu_1',\mu_2',\ldots,\mu_{N-1}')$$\in$$\Pi(\mathbf{M}_{fo})$ constructed through (\ref{MDP_controller}) as ${V}_{t,N}^{\pi'}(s^{t})$, respectively. For the base case, $t$$=$$N$$-$$1$, we have
        \begin{subequations}
    \begin{align}
       {V}_{N-1,N}^{\pi}(s^{N-1}) & = H^{\pi}(S_{N}|S^{N-1}=s^{N-1}) \label{eq:InductLastStep1-a} \\
        & = H^{\pi'}(S_{N}|S^{N-1}=s^{N-1}) \label{eq:InductLastStep1-b}
          \end{align}
\end{subequations}  
        by the definition of $\pi'$. In particular, the equality in (\ref{eq:InductLastStep1-b}) follows from the fact that we can construct an equivalent history-dependent controller on the fully observable MDP that achieves the same transition probabilities for any observation-based controller. Since the above equality holds for any $\pi$$\in$$\Pi(\mathbf{M})$, using the fact that  $\Pi(\mathbf{M})$$\subseteq$$\Pi(\mathbf{M}_{fo})$, we conclude
        \begin{subequations}
    \begin{align*}
        \sup_{\pi \in \Pi(\mathbf{M})} {V}_{N-1,N}^{\pi}(s^{N-1}) & \leq \sup_{\pi \in \Pi(\mathbf{M}_{fo})} {V}_{N-1,N}^{\pi}(s^{N-1}) . 
    \end{align*}
\end{subequations}

Now assume that, for any $s^{t+1}$$\in$$\mathcal{SH}^{t+1}$, we have
        \begin{subequations}
    \begin{align*}
        \sup_{\pi \in \Pi(\mathbf{M})} {V}_{t+1,N}^{\pi}(s^{t+1}) & \leq \sup_{\pi \in \Pi(\mathbf{M}_{fo})} {V}_{t+1,N}^{\pi}(s^{t+1}). 
    \end{align*}
\end{subequations}
Let $\overline{\pi}$$=$$(\overline{\mu}_1,\overline{\mu}_2,\ldots,\overline{\mu}_{N-1})$$\in$$\Pi(\mathbf{M}_{fo})$ be a history-dependent controller that attains the maximum on the right hand side of the above inequality. Then, for any $s^t$$\in$$\mathcal{SH}^t$, we have 
\begin{subequations}
    \begin{align}
        {V}_{t,N}^{\pi}(s^{t}) =& H^{\pi}(S_{t+1}|S^{t}=s^{t}) \nonumber \\
        & + \sum_{\substack{s^{t+1} \in \mathcal{S}\mathcal{H}^{t+1}}} Pr^{\pi}(s^{t+1}|s^{t}) {V}_{t+1,N}^{\pi}(s^{t+1}) \label{eq:InductMidStep1-a} \\
        \leq & H^{\pi}(S_{t+1}|S^{t}=s^{t}) \nonumber \\
        &  + \sum_{\substack{s^{t+1} \in \mathcal{S}\mathcal{H}^{t+1}}} Pr^{\pi}(s^{t+1}|s^{t}) {V}_{t+1,N}^{\overline{\pi}}(s^{t+1}) \label{eq:InductMidStep1-b} 
    \end{align}
\end{subequations}
due to Lemma \ref{writeValueFunction} and the induction hypothesis. We define a controller $\widetilde{\pi}$$=$$(\mu_1',\mu_2',\ldots,\mu_t',\overline{\mu}_{t+1},\ldots,\overline{\mu}_{N-1})$$\in$$\Pi(\mathbf{M}_{fo})$ which is a combination of the controllers $\pi'$ and $\overline{\pi}$. Note that
\begin{subequations}
    \begin{align}
        {V}_{t,N}^{\pi}(s^{t})
        \leq & H^{\widetilde{\pi}}(S_{t+1}|S^{t}=s^{t}) \nonumber \\
        &  + \sum_{\substack{s^{t+1} \in \mathcal{S}\mathcal{H}^{t+1}}} Pr^{\widetilde{\pi}}(s^{t+1}|s^{t}) {V}_{t+1,N}^{\widetilde{\pi}}(s^{t+1}) 
    \end{align}
\end{subequations}
since the controller $\widetilde{\pi}$ 
achieves the same transition probabilities with the controller $\pi$ for all state histories of length $t+1$ and the same value function for all state histories that has a length greater than $t+1$. Consequently, for any $s^t$$\in$$\mathcal{SH}^t$, we have
        \begin{subequations}
    \begin{align*}
        \sup_{\pi \in \Pi(\mathbf{M})} {V}_{t,N}^{\pi}(s^{t}) & \leq \sup_{\pi \in \Pi(\mathbf{M}_{fo})} {V}_{t,N}^{\pi}(s^{t}). \qquad \Box
    \end{align*}
\end{subequations}

\textbf{Proof of Proposition \ref{pmc_prop}:} The result follows from the fact that the controller $\mathbf{C}$ only allows deterministic memory transitions and that $\mathbf{D}_{\mathbf{M},k}[u_{\mathbf{C}}]$ is a Markov chain. By the definition of conditional entropy \cite[Chapter 2]{Cover}, 
\begin{flalign}
   & H^{\mathbf{C}}(S_t| S^{t-1})=\sum_{s^t\in \mathcal{S}\mathcal{H}^t} Pr^{\mathbf{C}}(s_t, s^{t-1})\log Pr^{\mathbf{C}}(s_t | s^{t-1}). && \raisetag{22pt}
\end{flalign}
Note that the summation on the right hand side of the above equation is over state histories. For any given state history $s^t$, there is a corresponding \textit{memory history} $(q_1,q_2,\ldots,q_t)$, where $q_k$$\in$$\mathcal{Q}$, for the controller $\mathbf{C}$. Let $\mathcal{MH}^t$ denote the set of all possible memory histories of length $t$$\in$$\mathbb{N}$. Then, by the law of total probability, we have
\begin{align*}
    Pr^{\mathbf{C}}(s_t | s^{t-1})&=\sum_{q^{t}\in \mathcal{MH}^t} Pr^{\mathbf{C}}(s_t,q_t | q^{t-1}, s^{t-1}) Pr^{\mathbf{C}}(q^{t-1}| s^{t-1}).
\end{align*}
Since the memory transitions are deterministic under $\mathbf{C}$$\in$$\mathcal{F}_k^{det}(\mathbf{M})$, by recursively expanding the right hand side of the above equality, it can be observed that $Pr^{\mathbf{C}}(q^{t-1}| s^{t-1})$$=$$1$ for a given state history realization $s^{t-1}$. Since for each state history realization $s^{t}$ on the POMDP $\mathbf{M}$ under the controller $\mathbf{C}$, there is a unique state history realization $(\langle s_1,q_1\rangle, \langle s_2,q_2\rangle, \ldots, \langle s_t,q_t\rangle)$ on the instantiation $\mathbf{D}_{\mathbf{M},k}[u_{\mathbf{C}}]$ of the induced pMC $\mathbf{D}_{\mathbf{M},k}$, we have 
\begin{align*}
    H^{\mathbf{C}}(S_{t} | S^{t-1})= H^{u_{\mathbf{C}}}(S_{\mathbf{M},k,t} | S_{\mathbf{M},k}^{t-1}).
\end{align*}
Finally, since the instantiation $\mathbf{D}_{\mathbf{M},k}[u_{\mathbf{C}}]$ constitutes an MC, as a result of the Markov property \cite{Cover}, we have $H^{u_{\mathbf{C}}}(S_{\mathbf{M},k,t} | S_{\mathbf{M},k}^{t-1})$$=$$H^{u_{\mathbf{C}}}(S_{\mathbf{M},k,t} | S_{\mathbf{M},k,t-1})$. $\Box$

\textbf{Proof of Lemma \ref{lemma:pMCmonotonic}:} We prove the claim by showing that, for any $k$$\in$$\mathbb{N}$, we have $E_{k-1,\max}$$\leq$$E_{k,\max}$.

Consider an arbitrary $(k$$-$$1)$-FSC $\mathbf{C}$$\in$$\overline{\mathcal{F}}_{k-1}(\mathbf{M})$ with the decision function $\gamma$. Let $k$-FSC $\mathbf{C}'$$\in$$\overline{\mathcal{F}}_{k}(\mathbf{M})$ be such that its decision function $\gamma'$ satisfies $\gamma'(a|q_{i},z)$$=$$\gamma(a|q_{i},z)$ for $i$$=$$1,\ldots,k$$-$$1$, and $\gamma'(a|q_{k},z)$$=$$\gamma(a|q_{k-1},z)$. Note that the state sequences in $\mathbf{M}$ under the controllers $\mathbf{C}$ and $\mathbf{C}'$ are the same. This is true since we can explicitly write down the decisions taken by the agent under the controllers $\mathbf{C}$ and $\mathbf{C}'$ thanks to the specific transition function given in \eqref{LastLoopDetFSC}. In particular, the sequence of decisions under the controller $\mathbf{C}$ is $(\gamma(a|q_1,z), \gamma(a|q_2,z), \ldots, \gamma(a|q_{k-1},z), \gamma(a|q_{k-1},z), \ldots)$, and the sequence of decisions under the controller $\mathbf{C}'$ is $(\gamma'(a|q_1,z), \gamma'(a|q_2,z), \ldots, \gamma'(a|q_{k},z), \gamma'(a|q_k,z), \ldots)$, which are the same by construction. Hence, the state sequences induced by these decision sequences are the same. Consequently, we have 
\begin{align*}
    \sum_{t=2}^{\infty}\beta^{t-2}H^{\mathbf{C}}(S_{t} | S^{t-1})=\sum_{t=2}^{\infty}\beta^{t-2}H^{\mathbf{C'}}(S_{t} | S^{t-1}).
\end{align*}
Since, for an arbitrary $(k$$-$$1)$-FSC $\mathbf{C}$$\in$$\overline{\mathcal{F}}_{k-1}(\mathbf{M})$, there exists a $k$-FSC $\mathbf{C}'$$\in$$\overline{\mathcal{F}}_{k}(\mathbf{M})$ that achieves the same entropy of state sequences in $\bf{M}$, we conclude that $E_{k-1,\max}$$\leq$$E_{k,\max}$. $\Box$

\section{} \label{section:Appendix2}
In this appendix, we describe a method to transform the finite horizon entropy maximization problem to infinite horizon entropy maximization problem with discount factor $\beta$$=$$1$. Even though the transformation requires us to use the discount factor $\beta$$=$$1$, the resulting POMDP includes a "sink state" which allows us to extend all theoretical results provided in the paper to finite horizon entropy maximization problem.

For a given POMDP $\mathbf{M}$ and a finite decision horizon $N$$\in$$\mathbb{N}$, we append the time as an additional state to the underlying transition system. In particular, instead of using $\mathcal{S}$, we use $(\mathcal{S}$$\times$$[N])$$\cup$$Sink$ as the state space, where $[N]$$=$$\{1,2,\ldots,N\}$. The initial state of the resulting POMDP is $s_I\times 1$, and the set of actions are the same as $\mathbf{M}$. The transition function $\overline{{P}}$ is 
\begin{align*}
    \overline{{P}}((s',t') | (s,t),a)=\begin{cases}
    {P}(s'|s,a) \ & \text{if} \ t'=t+1 \land t' < N\\
    1 \ & \text{if} \ (s',t')=Sink \land t= N\\
    1 \ & \text{if} \ (s',t')=(s,t)=Sink\\
    0 \ & \text{otherwise.}
    \end{cases}
\end{align*}
Intuitively, the process moves forward in time and get absorbed in the $Sink$ state after $N$-th stage. Now, on the resulting POMDP, all results provided for the infinite horizon setting can be extended to the finite horizon setting.

\vspace{-1cm}
\begin{IEEEbiography}[{\includegraphics[width=1in,height=1.25in,clip,keepaspectratio]{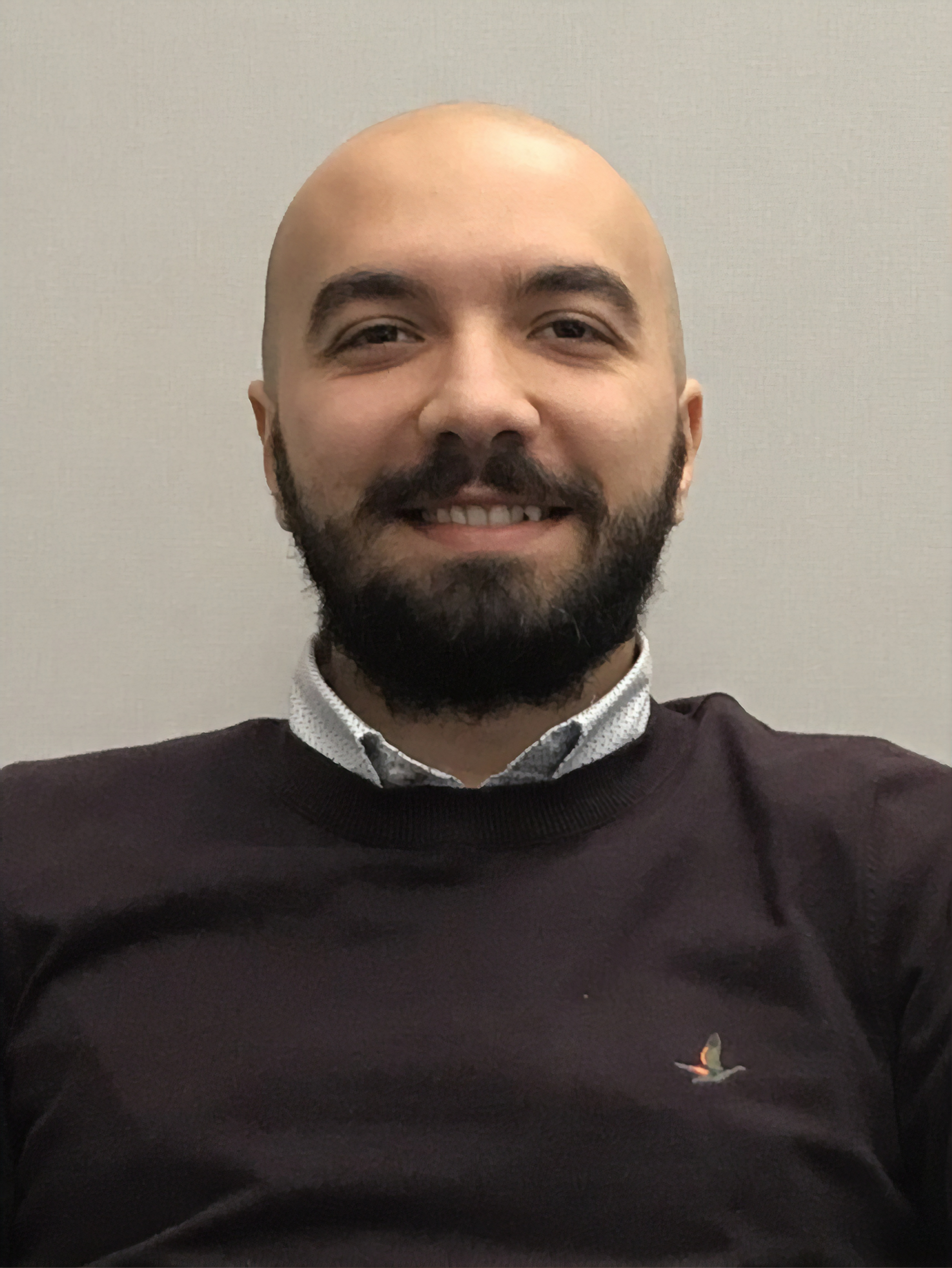}}]{Yagiz Savas} joined the Department of Aerospace Engineering at the University of Texas at Austin as a Ph.D. student in Fall 2017. He received his B.S. degree in Mechanical Engineering from Bogazici University in 2017. His research focuses on developing theory and algorithms that guarantee desirable behavior of autonomous systems operating in uncertain and adversarial environments.
\end{IEEEbiography}
\vspace{-1.1cm}
\begin{IEEEbiography}[{\includegraphics[width=1in,height=1.25in,clip,keepaspectratio]{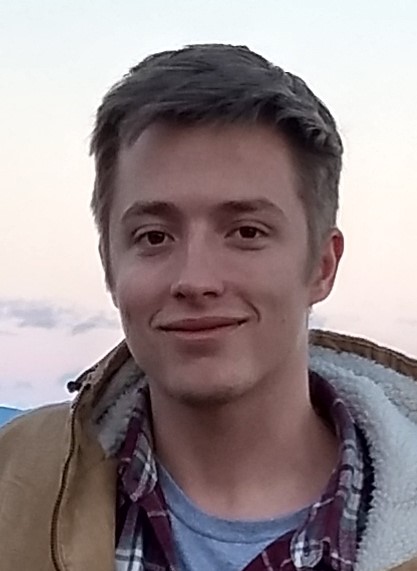}}]{Michael Hibbard} joined the Department of Aerospace Engineering at the University of Texas at Austin as a Ph.D. student in Fall 2018. He received his B.S. degree in Engineering Mechanics and Astronautics from the University of Wisconsin-Madison in Spring 2018. His research interests lie in the development of theory and algorithms providing formal guarantees for the mission success of autonomous agents acting in adversarial environments.
\end{IEEEbiography}
\vspace{-1.1cm}
\begin{IEEEbiography}[{\includegraphics[width=1in,height=1.25in,clip,keepaspectratio]{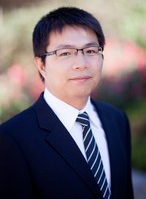}}]{Bo Wu}
received a B.S. degree from Harbin Institute of Technology, China, in 2008, an M.S. degree from Lund University, Sweden, in 2011 and a Ph.D. degree from the University of Notre Dame, USA, in 2018, all in electrical engineering. He is currently a postdoctoral researcher at the Oden Institute for Computational Engineering and Sciences at the University of Texas at Austin. His research interest is to apply formal methods, learning, and control in autonomous systems to provide privacy, security, and performance guarantees.
\end{IEEEbiography}
\vspace{-1.1cm}
\begin{IEEEbiography}[{\includegraphics[width=1in,height=1.25in,clip,keepaspectratio]{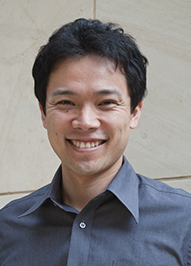}}]{Takashi Tanaka} received the B.S. degree from the University of Tokyo, Tokyo, Japan, in 2006, and the M.S. and Ph.D. degrees in Aerospace Engineering (automatic control) from the University of Illinois at Urbana-Champaign (UIUC), Champaign, IL, USA, in 2009 and 2012, respectively.
He was a Postdoctoral Associate with the Laboratory for Information and Decision Systems (LIDS) at the Massachusetts Institute of Technology (MIT), Cambridge, MA, USA, from 2012 to 2015, and a postdoctoral researcher at KTH Royal Institute of Technology, Stockholm, Sweden, from 2015 to 2017.
Currently, he is an Assistant Professor in the Department of Aerospace Engineering and Engineering Mechanics at the University of Texas at Austin.
His research interests include control theory and its applications; most recently the information-theoretic perspectives of optimal control problems.
\end{IEEEbiography}
\vspace{-1cm}
\begin{IEEEbiography}[{\includegraphics[width=1in,height=1.25in,clip,keepaspectratio]{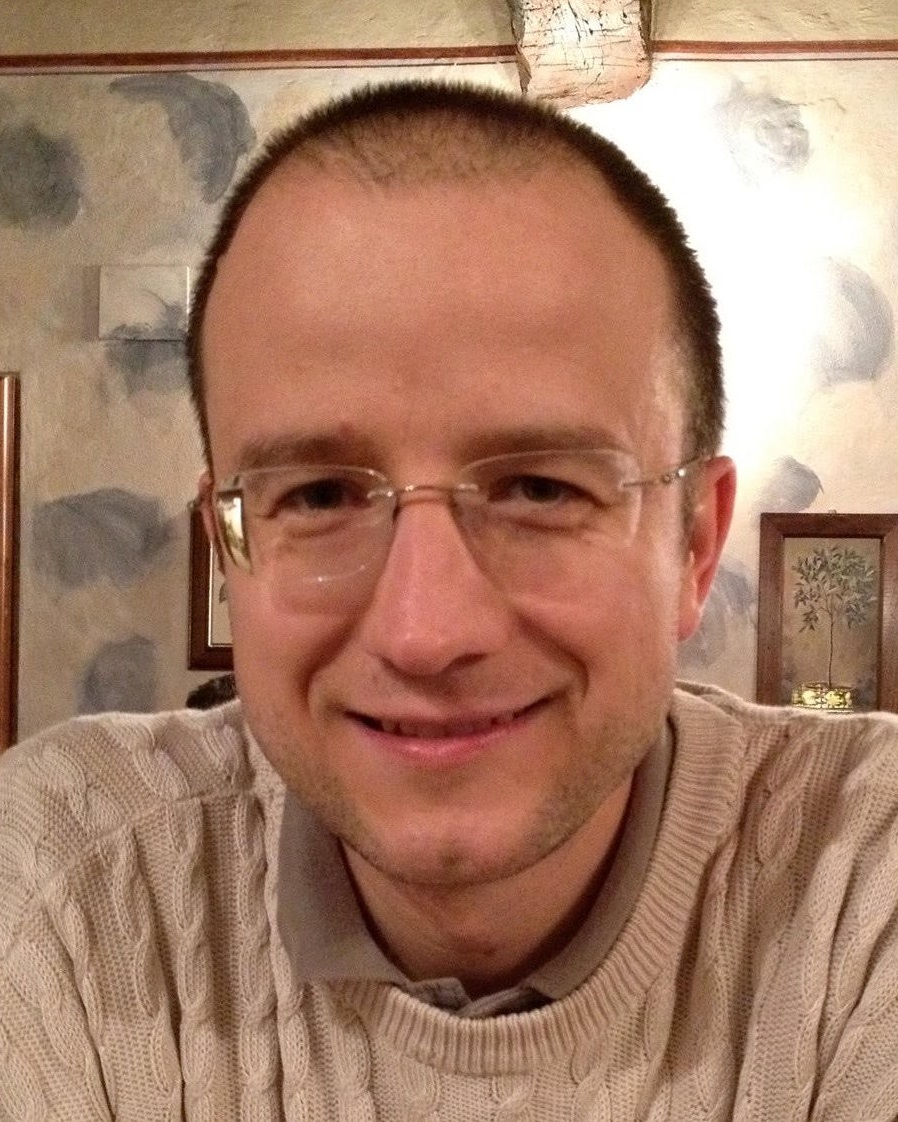}}]{Ufuk Topcu} joined the Department of Aerospace Engineering at the University of Texas at Austin as an assistant professor in Fall 2015. He received his Ph.D. degree from the University of California at Berkeley in 2008. He held research positions at the University of Pennsylvania and California Institute of Technology. His research focuses on the theoretical, algorithmic and computational aspects of design and verification of autonomous systems through novel connections between formal methods, learning theory and controls.
\end{IEEEbiography}

\end{document}